\documentclass[reqno]{amsart}

\usepackage{xcolor}

\allowdisplaybreaks
\newtheorem{theorem}{Theorem}[section]
\newtheorem{lemma}[theorem]{Lemma}
\newtheorem{corollary}[theorem]{Corollary}
\newtheorem{definition}[theorem]{Definition}
\newtheorem{proposition}[theorem]{Proposition}
\theoremstyle{definition}
\newtheorem{remark}[theorem]{Remark}

\numberwithin{equation}{section}

\numberwithin{equation}{section}

\newcommand{ \mr }{ \mathbb{R} }

\def\Xint#1{\mathchoice
    {\XXint\displaystyle\textstyle{#1}}%
     {\XXint\textstyle\scriptstyle{#1}}%
     {\XXint\scriptstyle\scriptscriptstyle{#1}}%
     {\XXint\scriptstyle\scriptscriptstyle{#1}}%
    \!\int}
\def\XXint#1#2#3{{\setbox0=\hbox{$#1{#2#3}{\int}$}
    \vcenter{\hbox{$#2#3$}}\kern-.5\wd0}}

\begin{document}

\title[Very weak solution in Orlicz space]{Existence of very weak solutions  to nonlinear elliptic equation with nonstandard growth and global weighted gradient estimates}

\author[S.-S. Byun, M. Lim]{Sun-Sig Byun\and Minkyu Lim}

\address{Sun-Sig Byun: Seoul National University, Department of Mathematical Sciences and Research Institute of Mathematics, Seoul 151-747, Korea}
\email{byun@snu.ac.kr}

\address{Minkyu Lim: Seoul National University,  Research Institute of Mathematics, Seoul 151-747, Korea}
\email{mk0314@snu.ac.kr}

\thanks{S.  Byun was supported by the National Research Foundation of Korea grant (NRF-2021R1A4A1027378). M. Lim was supported by the National Research Foundation of Korea grant (NRF-2022R1A2C1009312). }



\begin{abstract}
We study a general class of quasilinear elliptic equations with
nonstandard growth to prove the existence of a very weak solution
to such a problem. A key ingredient in the proof is a priori global
weighted gradient estimate of a very weak solution, where the right
hand side of the equation is the divergence of a vector-valued
function with low degree of integrability. To obtain this estimate,
we adopt a notion of reverse H\"older class of Muckenhoupt weights.
Another crucial part of the proof is a generalized weighted div-curl
lemma in the setting of Orlicz spaces.
\end{abstract}

\maketitle

\section{Introduction}
Throughout the paper, we mainly study the existence of a solution to the following equation on a bounded Lipschitz domain $\Omega \subset \mr^n \, (n \geq 2 )$
\begin{equation}\label{maineq3}
    \begin{cases}
        \mathrm{div\,}\mathcal{A}(x,Du)=\mathrm{div\,} \left( \frac{\varphi'(|\mathbf{f}|)}{|\mathbf{f}|} \mathbf{f} \right)  & \textrm{in}\ \Omega  \\
        u=0 & \textrm{on}\ \partial\Omega,
    \end{cases}
\end{equation}
when a given function $\mathbf{f} : \Omega \rightarrow \mr^n $ has a
 so low degree of integrability that we are below the natural exponent. Here, $\varphi \in C^1([0,\infty))
\cap C^2((0,\infty))$ with  $\varphi(0)=0$ is a so-called Young
function, which is convex and increasing (see \cite{RR}), and it
satisfies
\begin{equation}\label{youngcondi}
    \frac{1}{s_{\varphi}}\leq \frac{t
        \varphi''(t)}{\varphi'(t)}\leq s_{\varphi} \quad (t >0)
\end{equation}
for some constant $s_{\varphi}\geq 1$. Then, the natural assumptions
on the Carath\'eodory map $\mathcal{A} = \mathcal{A}(x,\xi):\Omega
\times \mr^n \rightarrow \mr^n $ in \eqref{maineq3} are the
following growth and monotonicity conditions:
\begin{equation}\label{moncondi}
    \begin{cases}
        |\mathcal{A}(x,\xi)| \leq L\varphi'(|\xi|)\\
        \langle \mathcal{A}(x,\xi)-\mathcal{A}(x,\zeta),\xi-\zeta\rangle \geq \nu
        \varphi''(|\xi|+|\zeta|)|\xi-\zeta|^2
    \end{cases}
\end{equation}
for a.e $x\in\Omega$ and all $\xi, \zeta \in\mr^n\backslash \{0\}$
with some constants $0 < \nu \leq L < + \infty$.

 In the case that
$\varphi(t)=t^{p} \, (1< p < \infty)$ and $
\mathcal{A}(x,\xi)=|\xi|^{p-2}\xi$, our equation \eqref{maineq3}
becomes a $p$-Laplacian type equation:
\begin{equation}\label{plaplace}
    \begin{cases}
        \mathrm{div\,} |Du|^{p-2} Du =\mathrm{div\,} |\mathbf{f}|^{p-2} \mathbf{f}   & \textrm{in}\ \Omega  \\
        u=0 & \textrm{on}\ \partial\Omega.
    \end{cases}
\end{equation}
Existence of a weak solution of the equation immediately follows
from the monotone operator theory if $\mathbf{f} \in L^{p}(\Omega)$,
see \cite[chapter 2]{Sh}. However, existence results for a distributional solution
to the equation \eqref{plaplace} have been obtained quite recently,
when the function $\mathbf{f}$ does not belong to $L^{p}(\Omega)$.
Especially, in \cite{BDS}, Bulicek, Diening and Schwarzacher proved
existence, uniqueness, and optimal regularity theory of a
distributional solution when the map $\mathcal{A}$ has linear growth
and satisfies asymptotically Uhlenbeck structure. In the later work
\cite{BS}, Bulicek and Schwarzacher showed an existence result for
general $p$-Laplacian type equations if $\mathbf{f} \in
L^{p-\varepsilon}(\Omega)$ for small $\varepsilon >0$. See also \cite{BBS, KL, MS}.  Our purpose
of this paper is to extend these result to a larger class of
quasilinear elliptic equations with nonstandard growth.

For the establishment of existence theory, one of the essential
requirements is to find proper weighted gradient estimates.
 Indeed, unweighted gradient estimates are
insufficient to establish solutions to \eqref{plaplace} when the
data $\mathbf{f}$ has such a low integrability. In a standard
approximation argument for the existence theory, unweighted $L^{q}$
estimates ($q < p$) might give a weak convergence of the sequence of
solutions in $L^{q}$ space, but these estimates alone do not ensure
that the limit of the sequence actually becomes a solution in
distributional sense. Although determining whether the limit become
a solution is usually performed via Minty method with the
monotonicity condition \eqref{moncondi}, this method is no longer
available, simply because of the fact that $q < p$. In the
aforementioned literature \cite{BDS, BS}, the authors overcome this
obstacle by replacing $L^{q}$ estimates with $L_{w}^{p}$ estimates,
where the weight $\omega$ is properly chosen in terms of
$\mathbf{f}$. The weighted gradient estimates recover the duality
relation between $Du $ and $\mathbf{f}$, which enables us to use the
condition \eqref{moncondi} to obtain the existence result.

 Here we are going to prove that there exists a
 positive constant $\delta_{*}$
depending on $n, s_{\varphi}, \nu$ and $L$ such that for some
constant $c$ independent of $u$ and $\mathbf{f}$, the estimate
\begin{equation}\label{mainesti}
    \int_{\Omega}    \varphi(|Du|)  \omega \, \mathrm{d}x   \leq c   \int_{\Omega}  \varphi(|\mathbf{f}|) \omega \, \mathrm{d}x
\end{equation}
holds for every Muckenhoupt weight $\omega $  in class $
A_{\frac{1}{1-\delta_{*}}} \cap RH_{1+\frac{1}{\delta_{*}}} $, see
\eqref{muckdef} and \eqref{rhdef} below for the classes $A_{p} \,
(1< p < \infty )$ and $RH_{s} (1< s < \infty )$. The
Muckenhoupt class $A_{p}$ is a set of weights for which the
Hardy-Littlewood maximal operator is bounded on $L_{\omega}^{p}$,
which is essentially used in the areas of harmonic analysis and
partial differential equations. Similarly, the reverse H\"older
class $RH_{s}$ consists of weight satisfying the reverse H\"older
inequality, which was first studied in \cite{CoF, G}. The relation
between these classes and weighted gradient estimates is profoundly
investigated in \cite{CF, CH, CN}.

As for the weighted gradient estimates \eqref{mainesti},
our results is new even in the special case that
$ \varphi(t) = t^{p}$. Indeed, for a solution
$u$ to the $p$-Laplacian type equation, the estimate
\eqref{mainesti} becomes
\begin{equation}\label{tempesti}
    \notag \int_{\Omega}  |Du|^{p} \omega \, \mathrm{d}x   \leq c   \int_{\Omega} |\mathbf{f}|^{p}  \omega \, \mathrm{d}x,
\end{equation}
where $\omega \in A_{\frac{1}{1-\delta_{*}}} \cap
RH_{1+\frac{1}{\delta_{*}}} $ with some positive constants
$\delta_{*}$ and $c$ independent of $u$ and $\mathbf{f}$. The above
weighted estimates for weak solutions to quasilinear elliptic and
parabolic equations have been intensively studied, with a BMO type
regularity assumption on the nonlinearity $\mathcal{A}$ with respect
to $x$ variables, as shown in \cite{BD, BP, MP, MP2, YCYY}. Recently, in
\cite{Ad} Adimurthi and Phuc proved a weighted estimate of the very
weak solution, which is a most natural concept of the solutions to
\eqref{plaplace} when the data $\mathbf{f}$ belongs to
$L_{\omega}^{q}(\Omega)$ with $q \geq p$, which covers the case that
$\mathbf{f} \not\in  L^{p}(\Omega) $. In that paper, this
restriction $q \geq p$ is unavoidable in the presence of the weight
$\omega$ as long as the choice of $\delta_{*}$ depends on
$A_{\infty}$ weight of $\omega $. Here we wish to relax this range
of $q $ in $ [p-\delta_{*}, \infty)$ by adapting the concept of a
reverse H\"older class and its related properties observed in
\cite{CN}. Consequently, we are able to find the most natural form
of the weighted gradient estimates for very weak solutions without
any regularity assumption on the nonlinearity $\mathcal{A}$.

Our paper is organized as follows. In the following
 section
we introduce the notion of very weak solutions, Orlicz-Sobolev
spaces, the Muckenhoupt class and reverse H\"older class to state
our main results.
 Section \ref{Sec3} contains various harmonic analysis tools,
 including
 Maximal function estimates, the properties of Muckenhoupt weights,
 an extrapolation lemma regarding weighted estimates,
 Calder\'on-Zygmund decomposition, and Lipschitz  truncation lemma.
In Section \ref{Sec4} we present comparison estimates of very weak
solutions to quasilinear equations and verify weighted gradient
estimates, Theorem \ref{mainthm2}. In the last section, we find a
variant of weighted div-curl lemma, adapted to our setting, to prove
the existence of a very weak solution of \eqref{maineq3} with the
desired estimate \eqref{mainesti2}

\section{Notation and main results}\label{Sec2}

 In this section, we present notations and preliminaries to describe our main results more precisely. For our quasilinear elliptic equation \eqref{maineq3}, we say that $u\in W_{0}^{1,1}(\Omega)$ is a very weak solution
to \eqref{maineq3} if $\varphi'(|Du|) \in  L^{1}(\Omega)$,
$\varphi'(\mathbf{|f|}) \in  L^{1}(\Omega)$ and
\begin{equation}\label{maineq2}
    \notag \int_{\Omega} \big\langle  \mathcal{A}(x,Dw), D\eta \big\rangle \,  \mathrm{d}x = \int_{\Omega} \big\langle  \frac{\varphi'(|\mathbf{f}|)}{|\mathbf{f}|} \mathbf{f} , D\eta \big\rangle  \,  \mathrm{d}x
\end{equation}
holds for every $ \eta \in C_{0}^{\infty}(\Omega)$.

 We then introduce standard notations which will be frequently throughout the paper.
 $B_\rho(y)$ denote the open ball
 with center $ y \in \mr^n$ and radius $\rho>0$. The
 intersection of the ball with  $\Omega\subset \mr^n$ will be denoted
 by $\Omega_\rho(y):=\Omega\cap B_\rho(y)$.   The complement of $\Omega$ is denoted by $\Omega^{c}:=\mr^n \backslash \Omega$.  We simply use
 the notation $B_\rho$ or $\Omega_\rho$ when the center $y$ is clear
 in the context. For an integrable
 function $v$ defined on a bounded measurable subset $E \subset
 \mathbb{R}^n$, we briefly denote to mean the integral average over $E$ by
 $$
 \overline{v}_E:=  \hspace{0.25em} \Xint-_E v(x)\; \mathrm{d}x =
 \frac{1}{|E|}\int_E v(x)\; \mathrm{d}x,
 $$
 where $|E|$ is the Lebesgue measure of $E$.   We extend $v$ by zero outside $F$ when $v$ is defined only on $F \subset E$. We denote the exponent conjugate to $p$ by $p':= \frac{p}{p-1}$ and denote the dual space of the vector space $X$ by $X^{*}$.

An auxiliary vector field is defined by
\begin{equation}\label{vtdef}
    \notag  V(\xi):=\left[ \frac{\varphi'(|\xi|)}{|\xi|} \right]^{\frac{1}{2}}\xi  \quad \left( \xi \in \mr^n\backslash \{0\} \right),
\end{equation}
which is a useful notion to describe the monotonicity property of
the vector field $\mathcal{A}(x,\xi)$ with respect to $\xi$ variables.
Indeed, if \eqref{youngcondi} is satisfied, then there holds
\begin{equation}\label{vtcomp}
    \notag \frac{1}{c}  \varphi''(|\xi|+|\zeta |)|\xi -\zeta |^2  \leq  |V(\xi )-V(\zeta )|^2   \leq c \varphi''(|\xi |+|\zeta |)|\xi -\zeta |^2  ,
\end{equation}
for some positive constant $c>1$, depending only on $n$ and $s_{\varphi}$, where  $ \xi  $ and $ \zeta $ belong to $  \mr^n\backslash \{0\} $. Then
\eqref{moncondi} yields
\begin{equation}\label{monot}
    \notag  \big\langle  \mathcal{A}(x,\xi )-\mathcal{A}(x,\zeta ) , \xi -\zeta  \big\rangle \geq c |V(\xi )-V(\zeta )|^2 \quad  \left( \xi  , \zeta  \in\mr^n\backslash \{0\} \right)
\end{equation}
for some positive constant $c$ depending on $n, \nu$ and
$s_{\varphi}$.

A domain under consideration is the following Lipschitz one regular
domain $\Omega$ throughout our paper.
\begin{definition} \label{lipschitz}
    We say that a domain $\Omega$ is a $(M, r_{0})$-Lipschitz domain if for every $x_{0} \in \partial \Omega$ and every $r \in (0, r_{0}]$, there exists a local coordinate  $y = (y_{1}, \cdots, y_{n}) \in B_{r}(0) $ with $y=0$
    at $x=x_{0}$,
    and a Lipschitz continuous function $\gamma : \mr^{n-1} \rightarrow \mr$ such that
    \begin{equation}\label{Lipdef1}
      \notag   \Omega \cap B_{r}(x_{0}) = \{ y = (y_{1}, \cdots , y_{n-1}, y_{n}) = (y', y_{n}) \in B_{r}(0) : y_{n} > \gamma (y')\}
    \end{equation}
    and
    \begin{equation}\label{Lipdef2}
      \notag   \sup_{|y'|<r, |z'|<r, y'\neq z'} \frac{|\gamma(y')-\gamma(z')|}{|y'-z'|} \leq M
    \end{equation}
    for some positive constant  $M \geq 1 $  and $r_{0}$.
\end{definition}
Note that for any $(M, r_{0})$-Lipschitz domain we have the following measure densitiy condition
\begin{eqnarray}\label{msrden}
      \min \left\{ {\frac{|B_{\rho}(x)\cap \Omega|}{|B_{\rho}(x)|}}, {\frac{|B_{\rho}(x)\cap \Omega^{c}|}{|B_{\rho}(x)|}}  \right\}  \geq \frac{1}{(n+1)!}\left(\frac{1}{M+1}\right)^{n-1}
\end{eqnarray}
for every ball $B_{\rho}(x)$ with $\rho\leq r_{0}$ and $x \in \partial \Omega$.

Next, we introduce the function spaces to which our solutions in the
main estimate \eqref{mainesti3} naturally belong. To do this, we
first recall weights and classes of weights commonly used. $ \omega
: \mr^n \rightarrow \mr $ is a weight if it is  a positive
measurable function and almost finite. A weight $\omega$ is an
$A_{p} $-weight $(1< p < \infty )$ if  $\omega \in
L^1_{\textrm{loc}}(\mr^n)$ and satisfies
\begin{eqnarray}\label{muckdef}
    [\omega]_{A_p} := \sup_{B_{\rho} \subset \mr^n} \left( \hspace{0.25em} \Xint-_{B_{\rho}} \omega(x)\; \mathrm{d}x \right) \left( \hspace{0.25em} \Xint-_{B_{\rho}} [\omega(x)]^{-\frac{1}{p-1}} \; \mathrm{d}x \right)^{p-1} < \infty.
\end{eqnarray}
Then we say that $\omega$ belongs to the Muckenhoupt class $A_{p}$
and the quantity $[\omega]_{A_p}$ is called to be the
$A_{p}$-constant of $\omega$. Similarly, the reverse H\"older class
$RH_{s} \, (1<s<\infty)$ is the collection of weights $\omega$ for
which
\begin{eqnarray}\label{rhdef}
    [\omega]_{RH_s} := \sup_{B_{\rho} \subset \mr^n} \left( \hspace{0.25em} \Xint-_{B_{\rho}} [\omega(x)]^{s}\; \mathrm{d}x \right)^{\frac{1}{s}} \left( \hspace{0.25em} \Xint-_{B_{\rho}}  \omega(x)  \; \mathrm{d}x \right)^{-1} < \infty.
\end{eqnarray}
As pointed out \cite{CN}, a function $\omega$ is in
$A_p$ for some $p > 1$ if and only if it is in $RH_s$ for some $s
> 1$. Then the reverse H\"older class provides another classification of
Muckenhoupt class. Meanwhile, Auscher and Martell has applied this
notion of the reverse H\"older class to a limited range of an
extrapolation in \cite{AM}, which is closely related to our
sub-natural gradient estimates, see for instance Lemma
\ref{extrapolation} below.

 We denote the set $L_{\omega}^{
    \varphi}(\Omega)$ for the collection of all measurable functions $v$ on
$\Omega$ for which
$$
\int_{\Omega} \varphi(|v(x)|) \omega(x) \; \mathrm{d}x < \infty.
$$
Indeed, if the function $\varphi$ satisfies condition \eqref{youngcondi}, $L_{\omega}^{
    \varphi}(\Omega)$ becomes a separable reflexive Banach space, where the norm is given
by
$$
\|v\|_{L_{ \omega}^{\varphi}(\Omega)}= \inf  \left\{ \lambda > 0 :
\int_{\Omega} \varphi \left( \frac{|v(x)|}{\lambda}\right) \omega(x) \; \mathrm{d}x \leq
1 \right\}.
$$
Then the function space $L_{\omega}^{\varphi}(\Omega)$ is a weighted
Orlicz space. See for instance [7]. $L_{\omega,
    \textrm{loc}}^{\varphi}(\Omega)$ is the union of the sets
$L_{\omega}^{\varphi}(\Omega' )$, where  $\Omega'  \subset\subset \Omega $. We simply write $L^{\varphi}(\Omega) := L_{\omega}^{\varphi}(\Omega)$ when $\omega = 1$ on $\mr^n$. We shortly remark that $(L^{\varphi}(\Omega))^{*} = L^{\varphi^{*}}(\Omega)$,  where $\varphi^{*}$ denotes the conjugate Young function of $\varphi$, that is,
$$
\varphi^{*}(t) :=  \sup_{s>0} \{ ts - \varphi(s) : s \in \mr \} \quad (t \geq 0).
$$
Moreover, if $q$ is sufficiently close to 1, $L^{{\varphi}^{q}}(\Omega)$ also becomes a separable
reflexive Banach space, where we have defined $\varphi^q(t) :=
[\varphi(t)]^{q}$ for $t\geq 0$, see Lemma \ref{orlicprop}.
Note that if $\omega \in A_{s_{\varphi}+1}  $, then we have $ L^{ \varphi }_{\omega}(\Omega)  \subset L^{ \varphi^{\frac{1}{s_{\varphi}+1}} }(\Omega) \subset L^{1 }(\Omega)  $.

 Then for any function $ u \in  L^{ \varphi }_{\omega}(\Omega)  \subset L^{1 }(\Omega)  $,
 we can consider its distributional derivatives. The Orlicz-Sobolev space
 $W_{\omega}^{1,\varphi}(\Omega)$ is a function space
consisting of all measurable functions $v\in
L_{\omega}^{\varphi}(\Omega)$ whose weak derivatives $Dv$ also belong to
$L_{\omega}^{\varphi}(\Omega;\mr^n)$.  The norm of
$W_{\omega}^{1,\varphi}(\Omega)$ is given by
$$
\|v\|_{W_{\omega}^{1,\varphi}(\Omega)}=
\|v\|_{L_{\omega}^{\varphi}(\Omega)}+\|Dv\|_{L_{\omega}^{\varphi}(\Omega)} .
$$
Note that
$\overline{C^{\infty}(B_{\rho})}=W_{\omega}^{1,\varphi}(B_{\rho})$ for
any ball $B_{\rho} \subset \mr^n$, where the completion is taken
with respect to the $W_{\omega}^{1,\varphi}(B_{\rho})$ norm.
$W_{\omega,0}^{1,\varphi}(\Omega)$ is defined as the closure of
$C_0^{\infty}(\Omega)$ in $W_{\omega}^{1,\varphi}(\Omega)$, where we have
denoted by $C_{0}^{\infty}(\Omega)$ to mean that the set of smooth
functions with compact support in $\Omega$.
We refer to \cite{DT} for a further discussion regarding this weighted
Orlicz-Sobolev spaces.

We now state our main theorems. One of our main goals is to prove
the following existence result for very weak solutions.
\begin{theorem}\label{mainthm}
    Assume \eqref{youngcondi} and  \eqref{moncondi}. Suppose that $\Omega $ is a $(M, r_{0})$-Lipschitz domain. Then there exists a
    small positive constant $\delta_{0}=\delta_{0}(n,s_{\varphi},  \nu,
    L, M)$ such that for all $\mathbf{f} \in L^{\varphi^{1-\delta_{0}}}(\Omega;\mr^n)$,
    there exists a very weak solution $u \in W_{0}^{1,\varphi^{1-\delta_{0}}} (\Omega)$ to
    \eqref{maineq3} with the estimates
    \begin{equation}\label{mainesti2}
        \int_{\Omega}  \varphi(|Du|)^{1-\delta_{0}} \, \mathrm{d}x   \leq c  \int_{\Omega}  \varphi(|\mathbf{f}|)^{1-\delta_{0}} \, \mathrm{d}x
    \end{equation}
    for some constant $c$ depending only on  $ n, s_{\varphi},
    \nu, L, M$ and $\frac{\mathrm{diam}(\Omega)}{r_{0}} $.

\end{theorem}

Another main result of this paper is an apriori weighted gradient
estimate in Orlicz setting, which plays a vital role in proving
Theorem \ref{mainthm}.

\begin{theorem}\label{mainthm2}
    Assume \eqref{youngcondi} and  \eqref{moncondi}. Suppose that $\Omega $ is a $(M, r_{0})$-Lipschitz domain.   Then there exists a
    small positive constant $\delta_{*}=\delta_{*}(n,s_{\varphi},  \nu,
    L, M)$ such that for all $\mathbf{f} \in L^{\varphi}_{\omega}(\Omega;\mr^n)$
    with $\omega \in A_{\frac{1}{1-\delta_{*}}} \cap RH_{1+\frac{1}{\delta_{*}}} $, any very weak solution $u \in W_{0}^{1,\varphi^{1-\delta_{*}}} (\Omega)$ of \eqref{maineq3}
    satisfies
    $$Du \in L^{ \varphi }_{\omega}(\Omega)$$
    with the estimate
    \begin{equation}\label{mainesti3}
      \notag  \int_{\Omega}  \varphi(|Du|)\omega(x) \, \mathrm{d}x   \leq c  \int_{\Omega}  \varphi(|\mathbf{f}|)\omega(x) \, \mathrm{d}x,
    \end{equation}
    where the constant $c$ depends only on  $ n, s_{\varphi},
    \nu, L, M, [\omega]_{A_{\frac{1}{1-\delta_{*}}}},  [\omega]_{RH_{1+\frac{1}{\delta_{*}}}}$ and $\frac{\mathrm{diam}(\Omega)}{r_{0}} $.

\end{theorem}

Although we establish this theorem as a tool for proving an existence result, the theorem itself has a significance
 in regularity theory of gradient estimates. For instance, one can directly apply this theorem
 along with Theorem \ref{extrapolation} in the next section,
   to prove the following  Calder\'on-Zygmund type estimate below the duality exponent.
\begin{corollary}\label{maincor}
    Assume \eqref{youngcondi} and  \eqref{moncondi}.  Suppose that $\Omega $ is a $(M, r_{0})$-Lipschitz domain. Then there exists a
    small positive constant $\delta_{*}=\delta_{*}(n,s_{\varphi}, \nu,
    L, M)$ such that for all $\mathbf{f} \in L^{\varphi^{q}}(\Omega;\mr^n)$
    with $q \in [1-\delta_{*}, 1+\delta_{*}]$, any very weak solution $u
    \in W_{0}^{1,\varphi^{1-\delta_{*}}}(\Omega)$ of \eqref{maineq3}
    satisfies
    $$ Du \in L^{\varphi^{q}}(\Omega)$$
    with the estimate
    \begin{equation}\label{mainesti3}
        \int_{\Omega}  [ \varphi(|Du|)]^q \, \mathrm{d}x   \leq c  \int_{\Omega} [ \varphi(|\mathbf{f}|)]^q \, \mathrm{d}x ,
    \end{equation}
    where the constant $c$ depends only on  $ n, s_{\varphi},  \nu, L, M$ and $\frac{diam(\Omega)}{r_{0}}$.
\end{corollary}
 Indeed, this unweighted estimates were proved   in the recent paper [4] under a certain capacitary thickness condition on the
 domain $\Omega$. Besides, if $q$ is in the range of $(1, \infty)$, then the estimate
 \eqref{mainesti3} becomes the standard Calder\'on-Zygmund type estimate,
 which holds under a suitable regularity condition on the map $x \mapsto A(x,z)$
 in \eqref{maineq3} and a sufficient flatness condition on the boundary $\partial \Omega$,
 as follows from  \cite{ Cho, Ve}. Collecting these facts, we can extend the exponent
 range to $q \in [1-\delta_{*}, \infty) $ in the above corollary by assuming some regularity
 requirements  on the nonlinearity and the boundary of the domain. A featured counterexample
 in [16] shows that  this extended exponent range is optimal range where Calder\'on-Zygmund type estimates  hold true.

\section{Preliminary Lemmas}\label{Sec3}
In this section, we list a few basic analysis tools to prove our main results.
We start with maximal function and Muckenhoupt weights. Let $\mathcal{M}$ be the Hardy-Littlewood maximal function defined by
\begin{eqnarray}\label{maxidef}
   \notag  \mathcal{M}(f)(y) = \sup_{\rho>0}  \hspace{0.25em} \Xint-_{B_{\rho}(y)}  |f(x)| \; \mathrm{d}x \quad ( x \in \mr^n )
\end{eqnarray}
for $f \in L^1_{ \textrm{loc}}(\mr^n)$. If $f\in L^1(\Omega)$ for some bounded domain $\Omega$, we extend $f$ by zero outside $\Omega$ so that the above definition makes sense. It is well known that if $\omega$ is an $A_{s}$-weight, then there exists a constant $c$ depending only on $s$ and $[\omega]_{A_s}$ such that
\begin{eqnarray}\label{muckdef2}
    \int_{B_{\rho}} [ \mathcal{M}(f)(x) ]^{s} \omega(x)  \, \mathrm{d}x  \leq c \int_{B_{\rho}} |f(x)|^{s} \omega(x)  \, \mathrm{d}x    \quad (B_{\rho} \subset \mr^n)
\end{eqnarray}
for every measurable function $f \in L^1_{ \textrm{loc}}(\mr^n)$ with the right hand side of \eqref{muckdef2} being finite \cite[Chapter 5]{St}.  We refer to \cite{BL} for the following lemma including standard classical theory regarding maximal functions and $A_{s}$ weights.
\begin{lemma}\label{maximucken}   Let $0<\tau<1$. Then for every nonnegative measurable function $f$ such that $f^{1-\delta} \in L^{1}(\mr^n)$ for some $\delta \in (0,\frac{1-\tau}{2}]$, there holds
    \begin{eqnarray}\label{maxiesti}
       \notag  \int_{B_{\rho}}  [ \mathcal{M} (f^{\tau}) (x)]^{\frac{1-\delta}{\tau}} \, \mathrm{d}x  \leq c  \int_{B_{\rho}}  [f(x)]^{1-\delta} \, \mathrm{d}x
    \end{eqnarray}
    for some positive constant $c$ depending on $n$ and $\tau$, whenever $B_{\rho} \subset \mr^n$. Moreover, $[\mathcal{M} (f^{\tau})]^{-\frac{\delta}{\tau}}$ is in the Muckenhoupt class $A_{\frac{1}{\tau}}$ with
    \begin{eqnarray}\label{muckenbound}
       \notag  [ [ \mathcal{M} (f^{\tau}) (x)]^{-\frac{\delta}{\tau}}]_{A_{\frac{1}{\tau}}} \leq c,
    \end{eqnarray}
    where the constant $c$ depends on $n$ and $\tau$.
\end{lemma}
Following so-called doubling property of Muckenhoupt weight will be useful.
\begin{lemma} \label{doublingweight} Let $\omega \in  A_{p} \cap RH_{s}$ for some $p>1$ and $s>1$. and let $E$ be a measurable subset of a ball $B$. Then
    \begin{eqnarray}\label{aprhproperty}
        \frac{1}{[\omega]_{A_{p}}}\left( \frac{|E|}{|B|} \right)^{p}  \leq \frac{\omega(E)}{\omega(B)} \leq [\omega]_{RH_{s}}\left( \frac{|E|}{|B|} \right)^{\frac{s-1}{s}}.
    \end{eqnarray}
\end{lemma}
\begin{proof}
The first inequality was proved in [5]. To prove the second inequality of \eqref{aprhproperty}, we observe
\begin{eqnarray}\label{doublingpf}
    \notag&&\hspace{-7mm} \left( \frac{|E|}{|B|}\right)^{\frac{s-1}{s}} = \left( \frac{\int_{B} \chi_{E} \,  \mathrm{d}x }{|B|}\right)^{\frac{s-1}{s}} = \left( \frac{\int_{B} \chi_{E}^{\frac{s}{s-1}} \,  \mathrm{d}x }{|B|}\right)^{\frac{s-1}{s}}   \\
    \notag&&  \geq  \frac{1}{|B|^{\frac{s-1}{s}}} \left(  \int_{B} \chi_{E} \omega  \,  \mathrm{d}x  \right)   \left(  \int_{B} \omega^{s} \,  \mathrm{d}x  \right) ^{-\frac{1}{s}} \\
    \notag&&  =   \frac{\omega(E)}{|B|^{\frac{s-1}{s}}} \cdot \frac{|B|}{\omega(B)}  \left(  \Xint-_{B}   \omega  \,  \mathrm{d}x  \right)   \left(  \Xint-_{B} \omega^{s} \,  \mathrm{d}x  \right) ^{-\frac{1}{s}}  \frac{1}{|B|^{\frac{1}{s}}}   \\
    \notag  &&  =   \frac{\omega(E)}{\omega(B)}  \left(  \Xint-_{B}  \omega  \,  \mathrm{d}x  \right)   \left(  \Xint-_{B} \omega^{s} \,  \mathrm{d}x  \right) ^{-\frac{1}{s}}  \geq   \frac{1}{[\omega]_{RH_{s}}}\frac{\omega(E)}{\omega(B)},
\end{eqnarray}
where we have used H\"older inequality for the inequality on the second line.
\end{proof}

Moreover, we have the following extrapolation lemma. A fundamental
form of extrapolation lemmas was discovered by Rubio de Francia in
\cite{R}, which is a very powerful tool in Harmonic Analysis. This
extrapolation lemma shows that the weighted norm inequality for one
single exponent propagates to the whole range $(1,\infty)$. Auscher
and Martell extended this result in \cite{AM}. They generalized Rubio
de Francia's lemma to the range $(p^{-}, p^{+})$ for any two
distinct positive constants $p^{-}$ and $p^{+}$ between $1$ and
$\infty$:
\begin{lemma} \cite[Corollary 4.10]{AM}\label{extrapolation}
    Given $1< p^{-} < p^{+}< \infty$, suppose that there exists $p_{0} \in (p^{-} , p^{+} )$ such that for every $\omega_{0} \in A_{\frac{p_{0}}{p^{-} }} \cap RH_{\left(\frac{p^{+}} {p_{0}}\right)'}$,
    \begin{equation}\label{weightassump}
     \notag    \int_{\Omega}  f(x)^{p_{0}} \omega_{0}(x) \, \mathrm{d}x   \leq c \int_{\Omega}  g(x)^{p_{0}} \omega_{0} (x)\, \mathrm{d}x.
    \end{equation}
    Then for every $p \in (p^{-} , p^{+} )$ and every $\omega \in A_{\frac{p}{p^{-} }} \cap RH_{\left(\frac{p^{+}} {p}\right)'}$ we have
    \begin{equation}\label{weightresult}
     \notag    \int_{\Omega}  f(x)^{p} \omega(x) \, \mathrm{d}x   \leq c \int_{\Omega}  g(x)^{p} \omega(x) \, \mathrm{d}x.
    \end{equation}
\end{lemma}
We remark that there are more general kinds of extrapolations. See
\cite{CF, CW} for the extrapolations in the weighted variable
exponent Lebesgue spaces and see \cite{CH} for the extrapolation in
the generalized Orlicz spaces.

The following basic properties of the Young function $\varphi$ are
found from \cite{ BL,Cho, HH}.
\begin{lemma} \label{orlicprop}
    Suppose $\varphi$ satisfies \eqref{youngcondi}. Then the following holds: \newline
    (a) For any $t>0$, we have
    \begin{equation}\label{orlicprop0}
       \notag  \frac{1}{s_{\varphi}}+1 \leq \frac{t \varphi'(t)}{\varphi(t)}\leq 1+s_{\varphi}.
    \end{equation} \newline
    (b) For any $a>0$, the function $\varphi_{a}(t):=\varphi(at)$ satisfies
    \begin{equation}\label{orlicprop1}
       \notag  \varphi_{a}(0)=0, \quad \frac{1}{s_{\varphi}}\leq \frac{t
            \varphi_{a}''(t)}{\varphi_{a}'(t)}\leq s_{\varphi} \quad (t >0).
    \end{equation} \newline
    (c) For any $0<\lambda \leq 1 $ and  $1 \leq \Lambda < \infty $, we have
    \begin{equation}\label{orlicprop2}
        \begin{cases}
            \lambda^{1+s_{\varphi}} \varphi(t)  \leq \varphi(\lambda t)\leq  \lambda^{(1/s_{\varphi})+1} \varphi(t)  \quad (t >0) \\
            \Lambda^{(1/s_{\varphi})+1} \varphi(t)  \leq \varphi(\Lambda t)\leq   \Lambda^{1+s_{\varphi}} \varphi(t)   \quad (t >0).
        \end{cases}
    \end{equation}
    (d) There exists a constant $\hat{\delta}=\hat{\delta} (s_{\varphi})>0$  such that for  all $q \in [1-\hat{\delta},1+\hat{\delta}]$,
    \begin{equation}\label{orlicprop3}
       \notag  \frac{1}{2s_{\varphi}} \leq \frac{t [\varphi^{q}(t)]''}{[\varphi^{q}(t)]'}\leq 2s_{\varphi} \quad (t >0) ,
    \end{equation}
    and $\varphi^{q}$ is increasing and convex.  \newline
    (e) For any
    $\varepsilon \in (0,1]$ and any $\delta \in [0, \hat{\delta}]$,
    there hold
    \begin{equation}\label{youngineq1}
      \notag   t [\varphi(s)]^{-\delta} \varphi'(s) \leq  \varepsilon [\varphi(t)]^{1-\delta} + c_{\varepsilon} [\varphi(s)]^{1-\delta}  \quad (t, s \geq 0)
    \end{equation}
    and
    \begin{equation}\label{youngineq2}
      \notag   t [\varphi(t)]^{-\delta} \varphi'(s)\leq  \varepsilon [\varphi(t)]^{1-\delta} + c_{\varepsilon} [\varphi(s)]^{1-\delta}  \quad (t, s \geq 0) ,
    \end{equation}
    where the constant $c_{\varepsilon}$ depends only on $s_{\varphi}$ and $\varepsilon$.
\end{lemma}

We next introduce an weighted Orlicz-Sobolev-Poincar\'e type inequality. We generalized the unweighted inequality found in \cite[Theorem 7]{DE}.
\begin{lemma}\label{poincare} Suppose $\varphi$ satisfies \eqref{youngcondi} and let $\omega \in A_{1+\frac{1}{s_{\varphi}}}$. Then there exists a constant $\theta=\theta(n, s_{\varphi}) > 1 $  such that for any $v \in W^{1, \varphi }(B_\rho)$, there holds
    \begin{eqnarray}\label{poin1}
        \left( \frac{1}{\omega(B_{\rho})} \int_{B_{\rho}} \left[ \varphi \left(\frac{|v-v_{\rho}|}{\rho}\right)   \right]^{\theta} \omega \, \mathrm{d}x \right)^{\frac{1}{\theta}}  \leq    \frac{c}{\omega(B_{\rho})}  \int_{B_{\rho}}   \varphi(|Dv|)   \omega \, \mathrm{d}x
    \end{eqnarray}
    for some positive constant $c$ depending only on $n, s_{\varphi}$ and $[\omega]_{A_{1+\frac{1}{s_{\varphi}}}}$.
\end{lemma}
\begin{proof}
    For almost every $x \in B_{\rho}$, consider the following classical formula
    \begin{eqnarray}\label{poinpf1}
        |v(x)-v_{\rho}| \leq c \int_{B_{\rho}}  \frac{|Dv(y)|}{|x-y|^{n-1}}    \, \mathrm{d}y
    \end{eqnarray}
     presented in \cite[Lemma 7.16]{GT}, where
    the constant $c$ depends only on $n$. Let $\varepsilon \in (0,1]$ and denote
     $$
     D_{i}(x) := B_{2^{-i}\varepsilon \rho}(x) \backslash B_{2^{-i-1}\varepsilon \rho}(x)
     $$
     for integer $i \geq 0$. Then \eqref{poinpf1} leads to
     \begin{eqnarray}\label{poinpf2}
        \notag &&  \hspace{-10mm} |v(x)-v_{\rho}| \leq    c\int_{ B_{\varepsilon \rho}(x) }  \frac{|Dv(y)|}{|x-y|^{n-1}}    \, \mathrm{d}y   + c\int_{B_{\rho}  \backslash B_{\varepsilon \rho}(x) }  \frac{|Dv(y)|}{|x-y|^{n-1}}    \, \mathrm{d}y   \\
        \notag &&   \hspace{-5mm} \leq    c \left[ \sum_{i=0}^{\infty}
        \left( \frac{2^{i}}{\varepsilon \rho} \right)^{n-1} \int_{ D_{i}(x) }  |Dv(y)|  \, \mathrm{d}y   +  \frac{1}{(\varepsilon \rho)^{n-1}} \int_{B_{2\rho}(x)    }  |Dv(y)|     \, \mathrm{d}y \right]  \\
        \notag &&   \hspace{-5mm} \leq    c \rho \left[ \sum_{i=0}^{\infty}
        \frac{\varepsilon}{2^{i}}   \hspace{0.25em} \Xint-_{ B_{2^{-i}\varepsilon \rho}(x) }  |Dv(y)|  \, \mathrm{d}y   +  \frac{1}{\varepsilon ^{n-1}} \hspace{0.25em} \Xint-_{B_{2\rho}(x) }  |Dv(y)|    \, \mathrm{d}y \right]\\
        &&   \hspace{-5mm} \leq    c \varepsilon \rho  \mathcal{M} (|Dv|) (x) + \frac{c \rho }{\varepsilon ^{n-1}} \hspace{0.25em} \Xint-_{B_{2\rho}(x) }  |Dv(y)|   \, \mathrm{d}y
     \end{eqnarray}
     for some constant $c$ depending only on $n$. As remarked at the beginning of  this section, $|Dv|$ was extended by zero outside $B_{\rho}$. Since $\varphi$ is increasing, with the property \eqref{orlicprop2} and the above inequality \eqref{poinpf2} we obtain
     \begin{align}\label{poinpf3}
       \varphi\left(\frac{|v(x)-v_{\rho}|}{\rho} \right) \leq  c \varepsilon^{\frac{1}{s_{\varphi}}} \varphi \Big(    \mathcal{M} (|Dv|) (x)  \Big) +  \frac{c}{\varepsilon^{s_{\varphi}(n-1)}} \varphi\left(     \hspace{0.25em} \Xint-_{B_{2\rho}(x) }  |Dv(y)|   \, \mathrm{d}y \right).
     \end{align}
     Now, if we choose
     \begin{equation}\label{poinpf4}
      \notag  \theta:= 1+\frac{ 1}{ s_{\varphi}^2(n-1)}  \quad \textrm{and} \quad \varepsilon := \left[  \frac{\varphi\left(     \hspace{0.25em} \Xint-_{B_{2\rho}(x) }  |Dv(y)|   \, \mathrm{d}y \right)}{ \varphi \Big(    \mathcal{M} (|Dv|) (x)  \Big)} \right]^{\frac{s_{\varphi}}{1+ s_{\varphi}^2(n-1)}    }  \leq 1,
     \end{equation}
     then the inequality \eqref{poinpf3} becomes
     \begin{eqnarray}\label{poinpf5}
              \left[ \varphi\left(\frac{|v(x)-v_{\rho}|}{\rho} \right) \right]^{\theta}    \leq    c  \varphi \Big(    \mathcal{M} (|Dv|) (x)  \Big)   \left[\varphi\left(     \hspace{0.25em} \Xint-_{B_{2\rho}(x) }  |Dv(y)|   \, \mathrm{d}y \right)\right]^{\theta-1}
     \end{eqnarray}
     for a.e. $x \in B_{\rho}$. Using the H\"older inequality and the fact that $\omega \in A_{1+\frac{1}{s_{\varphi}}}$, we get
     \begin{eqnarray}\label{poinpf6}
        \notag && \hspace{-7mm}\varphi\left(     \hspace{0.25em} \Xint-_{B_{2\rho}(x) }  |Dv(y)|   \, \mathrm{d}y \right)   \\
        \notag &&   \leq     \varphi\left(  \left( \hspace{0.25em} \Xint-_{B_{2\rho}(x) }  |Dv(y)|^{1+\frac{1}{s_{\varphi}}}  \omega \, \mathrm{d}y \right)^{\frac{s_{\varphi}}{s_{\varphi}+1}}   \left(     \hspace{0.25em} \Xint-_{B_{2\rho}(x) }  \omega^{-s_{\varphi}}  \, \mathrm{d}y  \right)^{\frac{1}{s_{\varphi}+1}}  \right)  \\
        \notag&&   \leq     \varphi\left(         \left( \hspace{0.25em} \Xint-_{B_{2\rho}(x) }  |Dv(y)|^{1+\frac{1}{s_{\varphi}}}  \omega \, \mathrm{d}y \right)^{\frac{s_{\varphi}}{s_{\varphi}+1}} \left(      \frac{[w]_{A_{1+ \frac{1}{s_{\varphi}}}} | B_{2\rho}(x) | }{\omega(B_{2\rho}(x))}\right)^{\frac{s_{\varphi}}{s_{\varphi}+1}}  \right)\\
        \notag&&   =     \varphi\left(         \left( \frac{  [w]_{A_{1+ \frac{1}{s_{\varphi}}}}}{\omega(B_{2\rho}(x))} \int_{B_{2\rho}(x) }  |Dv(y)|^{1+\frac{1}{s_{\varphi}}}  \omega \, \mathrm{d}y \right)^{\frac{s_{\varphi}}{s_{\varphi}+1}}            \right) \\
        &&   \leq             \frac{ c }{\omega(B_{2\rho}(x))} \int_{B_{2\rho}(x) }  \varphi\left( |Dv(y)| \right)  \omega \, \mathrm{d}y
     \end{eqnarray}
     for some constant $c$ depending only on $n, s_{\varphi}$ and $ [w]_{A_{1+ \frac{1}{s_{\varphi}}}}$, where we have used the Jensen inequality with a convex function $\varphi(t^{\frac{s_{\varphi}}{s_{\varphi}+1}})$ for the last inequality of \eqref{poinpf6}.  Meanwhile, we have the following maximal function estimates
     \begin{eqnarray}\label{poinpf7}
        \hspace{-3mm}      \int_{\mr^n}  \varphi \Big(    \mathcal{M} (|Dv|) (x)  \Big)     \omega \, \mathrm{d}x      \leq    c \int_{ \mr^n }  \varphi\left( |Dv(x)| \right)  \omega \, \mathrm{d}x
     \end{eqnarray}
     for any weight $\omega \in A_{1+\frac{1}{s_{\varphi}}}$  as follows from \cite[Theorem 2.1.1]{KK}.

     Combining \eqref{poinpf5}, \eqref{poinpf6} and \eqref{poinpf7}, we obtain
     \begin{eqnarray}\label{poinpf8}
        \notag && \hspace{-3mm}    \frac{1}{\omega(B_{\rho})} \int_{B_{\rho}} \left[ \varphi \left(\frac{|v(x)-v_{\rho}|}{\rho}\right)   \right]^{\theta} \omega(x) \, \mathrm{d}x    \\
        \notag && \hspace{-6mm}  \leq       \frac{c}{\omega(B_{\rho})} \int_{B_{\rho}} \varphi \Big(    \mathcal{M} (|Dv|) (x)  \Big)  \omega(x)    \left(  \frac{ 1 }{\omega(B_{2\rho}(x))} \int_{B_{2\rho}(x) }  \varphi\left( |Dv(y)| \right)  \omega(y) \, \mathrm{d}y  \right)^{\theta-1}   \, \mathrm{d}x  \\
        \notag&&  \hspace{-6mm}   \leq       \frac{c}{\omega(B_{\rho})} \int_{B_{\rho}} \varphi \Big(    \mathcal{M} (|Dv|) (x)  \Big)  \omega(x)   \, \mathrm{d}x   \left(  \frac{ 1 }{\omega(B_{\rho})} \int_{B_{\rho}  }  \varphi\left( |Dv(y)| \right)  \omega(y) \, \mathrm{d}y  \right)^{\theta-1}   \\
        \notag &&  \hspace{-6mm}  \leq    \left( \frac{  c }{\omega(B_{\rho})} \int_{ B_{\rho} }  \varphi\left( |Dv(x)| \right)  \omega(x) \, \mathrm{d}x \right)^{\theta}
     \end{eqnarray}
     for some constant $c$ depending only on $n, s_{\varphi}$ and $ [w]_{A_{1+ \frac{1}{s_{\varphi}}}}$. Then the desired estimate \eqref{poin1} follows.
\end{proof}

The following  Lipschitz truncation lemma plays a crucial role in proving div-curl lemma \ref{divcurl}.
\begin{lemma}
    \label{liptrun}
    Suppose $\Omega \subset \mr^n $ is a $(M, r_{0})$-Lipschitz domain. For $v  \in W_{0}^{1,1}(\Omega )$, we write
    $$
    \quad E_{\lambda}:=  \{x \in \Omega :   \mathcal{M}(Dv)(x) >\lambda
    \}  \quad (\lambda>0).
    $$
    Then there exist a Lipschitz function $v_\lambda \in W_{0}^{1,\infty} (\Omega) $
    and a positive constant $c$ depending on $n, s_{\varphi}, M$ and $\frac{\mathrm{diam}(\Omega)}{r_{0}}$
    such that
    $$v_\lambda(x)=v(x),   \quad  Dv_\lambda(x)=Dv(x) \quad    \quad \mathrm{a.e.} \ x \in \Omega \backslash E_{\lambda}  $$
    and that the estimate
    $ |Dv_{\lambda}(x) |  \leq c\lambda$
    holds for a.e. $ x \in \Omega $. Furthermore, if $\varphi$ satisfies
    \eqref{youngcondi} and $Dv  \in L_{\omega}^{\varphi}(\Omega )$ with
    $\omega \in A_{1+\frac{1}{s_{\varphi}}} $, then we have
    \begin{eqnarray}
        &&
         \int_{\Omega} \varphi ( |D v_\lambda| )  \omega(x) \, \mathrm{d}x   \leq c \int_{\Omega}  \varphi (   |D v|  )    \omega(x) \, \mathrm{d}x \label{lipsestidif} , \label{lipsesti} \\
        &&  \hspace{-5mm} \int_{\Omega} \varphi ( |D v - D v_\lambda| )  \omega(x) \, \mathrm{d}x   \leq c \int_{E_{\lambda}}  \varphi (   |D v|  )    \omega(x) \, \mathrm{d}x \label{lipsestidif}
    \end{eqnarray}
    for some positive constant $c$ depending only on $n, s_{\varphi}, M, \frac{\mathrm{diam}(\Omega)}{r_{0}}$ and $[\omega]_{A_{1+\frac{1}{s_{\varphi}}}}$.
\end{lemma}
\begin{proof}
         We define the set
        $$
        \quad F_{\lambda}:=  E_{\lambda}^{c} = \{x \in \Omega : \mathcal{M}(Dv)(x) \leq \lambda \} \cup \Omega^{c}
        $$
        for every $\lambda>0$, which is closed by the lower semicontinuity of the function $\mathcal{M}(Dv)$ in $\mr^n$. Without loss of generality, we assume that $E_{\lambda} \neq \emptyset $ and let $w \in W^{1,1}(\mr^n)$ denote the extension of $v$ by zero outside $\Omega$.  For every balls $B_{R}(x_{0}) \subset \mr^n$ with $R>0$ and $x_{0} \in \mr^n$, we denote $R_{k}:= 2^{1-k}R$. Then for every Lebesgue point $\xi \in B_{R}(x_{0})$ of $w$, we have
        \begin{eqnarray}\label{ptesti1}
            \notag &&  \hspace{-5mm}  |w(\xi) -  \bar{w}_{B_{R}(x_{0})}|  \leq  \sum_{k=0}^{\infty}   \left|  \bar{w}_{B_{R_{k+1}}(\xi) } - \bar{w}_{B_{R_{k}}(\xi) }  \right| +  \left|  \bar{w}_{B_{2R}(\xi)} -  \bar{w}_{B_{R}(x_{0})} \right|  \\
            \notag &&    \hspace{12mm}  \leq  \sum_{k=0}^{\infty}  \left[  \hspace{0.25em} \Xint-_{B_{R_{k+1}}(\xi)} \left|  w   - \bar{w}_{B_{R_{k}}(\xi) } \right|  \, \mathrm{d}x \right]   +   \hspace{0.25em} \Xint-_{B_{R}(x_{0})} \left|  w   - \bar{w}_{B_{2R}(\xi) } \right|  \, \mathrm{d}x \\
            \notag &&    \hspace{12mm}  \leq  \sum_{k=0}^{\infty} 2^n \left[  \hspace{0.25em} \Xint-_{B_{R_{k}}(\xi)} \left|  w   - \bar{w}_{B_{R_{k}}(\xi) } \right|  \, \mathrm{d}x \right]   +   2^n \hspace{0.25em} \Xint-_{B_{2R}(\xi)} \left|  w   - \bar{w}_{B_{2R}(\xi) } \right|  \, \mathrm{d}x \\
            \notag &&    \hspace{12mm} \leq  c \left[  \sum_{k=0}^{\infty} 2^{-k} R  \left(  \hspace{0.25em} \Xint-_{B_{R_{k}}(\xi)} \left|   Dw \right|  \, \mathrm{d}x \right)   +   R  \hspace{0.25em} \Xint-_{B_{2R}(\xi)} \left|  Dw  \right|  \, \mathrm{d}x \right]  \\
            &&    \hspace{12mm}  \leq c R \mathcal{M}(Dv)(\xi),
        \end{eqnarray}
        for some constant $c$ depending on $n$, where we have used the Poincar\'e inequality for the fourth inequality. Then for every Lebesgue point $\xi, \zeta \in  F_{\lambda} \cap \Omega$ of $w$, we obtain
        \begin{equation}\label{ptesti2}
          \notag      \frac{ |w(\xi) -  w(\zeta)|}{|\xi-\zeta|}  \leq c    [ \mathcal{M}(Dv)(\xi) +  \mathcal{M}(Dv)(\zeta) ] \leq c
               \lambda,
        \end{equation}
        where we have taken $x_{0}:= \frac{\xi + \zeta}{2}$ and $R=|\xi-\zeta|$ in \eqref{ptesti1}. In the case when $\xi \in  F_{\lambda} \cap \Omega$ but
        $\zeta \in \Omega^{c}$, we further assume that $|\xi- \zeta |\geq
        r_{0}$. Since $w \in W_{0}^{1,1}(B_{2 \,
            \mathrm{diam}(\Omega)}(\xi))$, we have
        \begin{eqnarray}\label{ptesti3}
            \notag &&  \hspace{-5mm}   \frac{ |w(\xi) -  w(\zeta)|}{|\xi- \zeta|}   \leq \frac{ |w(\xi) |}{r_{0}}  \leq  \frac{|w(\xi) -  \bar{w}_{B_{2 \, \mathrm{diam}(\Omega)}(\xi)}|}{r_{0}} + \frac{|\bar{w}_{B_{2 \, \mathrm{diam}(\Omega)}(\xi)}|}{r_{0}}   \\
            \notag &&  \underset{\eqref{ptesti1}}{\leq} \frac{ c \, \mathrm{diam}(\Omega) }{r_{0}}  \left[  \mathcal{M}(Dv)(\xi)  +     \hspace{0.25em} \Xint-_{B_{2 \, \mathrm{diam}(\Omega)}(\xi)}|Dw|  \, \mathrm{d}x    \right]  \\
            \notag  && \hspace{2mm} \leq  c  \mathcal{M}(Dv)(\xi) \leq c \lambda,
        \end{eqnarray}
        for some constant $c$ depending on $n$ and $\frac{\mathrm{diam}(\Omega)}{r_{0}} $. On the other hand,
        if $|\xi- \zeta | <  r_{0}$, then there exists a point $\zeta_{0} \in \partial \Omega $
        such that $|\xi- \zeta_{0} | <  r_{0}$. For every Lebesgue point $\xi$ of $w$,
        the Poincar\'e inequality with measure density condition \eqref{msrden} gives
        \begin{eqnarray}\label{ptesti4}
            \notag &&  \hspace{-5mm}    \frac{ |w(\xi) -  w(\zeta)|}{|\xi- \zeta|}   \leq  \frac{ |w(\xi) |}{|\xi- \zeta_{0}|}  \leq  \frac{ |w(\xi) - \bar{w}_{B_{|\xi- \zeta_{0}|}(\xi)} |}{|\xi- \zeta_{0}|}   + \frac{ |\bar{w}_{B_{|\xi- \zeta_{0}|}(\xi)}| }{|\xi- \zeta_{0}|}    \\
            \notag &&  \leq c      \mathcal{M}(Dv)(\xi)     +      \hspace{0.25em} \Xint-_{    \Omega^{c}    \cap B_{|\xi- \zeta_{0}|}(\xi) }   \frac{|w- \bar{w}_{B_{|\xi- \zeta_{0}|}(\xi)}|}{|\xi- \zeta_{0}|}    \, \mathrm{d}x      \\
            \notag &&   \leq c      \mathcal{M}(Dv)(\xi)   \leq c \lambda
        \end{eqnarray}
        for some constant $c$ depending on $n, s_{\varphi}$ and $M$. Then $w$ becomes a $c^{*}
         \lambda $-Lipschitz function on the closed set
        $F_{\lambda}$ for some constant $c^{*} $ depending on $n,
        s_{\varphi}, M $ and $\frac{\mathrm{diam}(\Omega)}{r_{0}}$.
        Therefore, there is a $c^{*}  \lambda $-Lipschitz
        extension $v_{\lambda} \in W^{1,\infty}(\mr^n)$ of $w$ by Kirszbraun
        extension theorem, see for instance \cite[Chapter 3]{EG}. Since
        $v_{\lambda}$ vanishes outside $\Omega$,  $v_\lambda \in W_{0}^{1,\infty} (\Omega) $.

        Now we check that the inequalities \eqref{lipsesti} and \eqref{lipsestidif} are valid. Observe that
        \begin{eqnarray}\label{lipesti1}
            \notag &&  \hspace{-5mm}  \int_{\Omega} \varphi ( |  D v_\lambda| )  \omega(x) \, \mathrm{d}x    \leq c  \int_{\Omega \backslash E_{\lambda}} \varphi ( |  D v | )  \omega(x) \, \mathrm{d}x + c    \varphi (\lambda ) w(E_{\lambda})  \\
            \notag &&  \leq c  \int_{\Omega \backslash E_{\lambda}} \varphi ( |  D v | )  \omega(x) \, \mathrm{d}x +  c\varphi (\lambda ) \int_{ E_{\lambda} } \frac{\varphi ( \mathcal{M}(Dv) )}{\varphi (\lambda )}  \omega(x) \, \mathrm{d}x    \\
            \notag &&  \hspace{-2mm}  \underset{\eqref{poinpf7}}{\leq}   c \int_{\Omega}  \varphi (   |D v|  )    \omega(x) \, \mathrm{d}x
        \end{eqnarray}
    for some constant $c$ depending on $n, s_{\varphi}, M, \frac{\mathrm{diam}(\Omega)}{r_{0}}$ and $[\omega]_{A_{1+\frac{1}{s_{\varphi}}}}$. This proves the estimate \eqref{lipsesti}. Similarly, since $Dv_{\lambda} = Dv$ on a.e. $x \in \Omega \backslash E_{\lambda}$, we have
    \begin{eqnarray}\label{lipesti2}
        \notag &&  \hspace{-5mm}  \int_{\Omega} \varphi ( |D v - D v_\lambda| )  \omega(x) \, \mathrm{d}x    = \int_{E_{\lambda}} \varphi ( |D v - D v_\lambda| )  \omega(x) \, \mathrm{d}x  \\
        \notag &&  \leq c  \int_{E_{\lambda}} \varphi ( |D v | )  \omega(x) \, \mathrm{d}x  + c \int_{E_{\lambda}} \varphi ( |  D v_\lambda| )  \omega(x) \, \mathrm{d}x    \\
        &&  \leq c  \int_{E_{\lambda}} \varphi ( |D v | )  \omega(x) \, \mathrm{d}x  +     c    \varphi (\lambda ) w(E_{\lambda})
    \end{eqnarray}
for some constant $c$ depending on $n, s_{\varphi}, M$ and $
\frac{\mathrm{diam}(\Omega)}{r_{0}}$. To estimate the second term of
\eqref{lipesti2}, we choose $r_{x}$ such that
$$
\lambda <\hspace{0.25mm} \Xint-_{B_{r_{x}}(x)}  |D v |  \,
\mathrm{d}x \leq 2\lambda
$$
for every point $x $ in $ E_{\lambda} $. Then the collection
$\{B_{r_{x}}(x)\}$ becomes a Besicovich cover of the set
$E_{\lambda}$.
 By Besicovich covering theorem, there is a countable subcover $\{B_{i} : i \in \mathbb{N} \}$ of the collection $\{B_{r_{x}}(x)\}$ which still covers $E_{\lambda}$, where each point $x$ of $E_{\lambda}$ is contained in at most $c(n)$ balls of $\{B_{i}\}$ with the constant $c(n)$ depending only on $n$. Then we have
 \begin{eqnarray}\label{lipesti3}
    \notag &&  \hspace{-5mm} \varphi (\lambda ) w(E_{\lambda})  \leq   \sum_{i \in \mathbb{N} }  \varphi (\lambda ) w(B_{i})  \leq  \sum_{i \in \mathbb{N} }   \varphi \left( \hspace{0.25mm} \Xint-_{B_{i}}  |D v |  \, \mathrm{d}x  \right) w(B_{i}) \\
    \notag &&  \leq  \sum_{i \in \mathbb{N} }  \varphi\left(  \left( \hspace{0.25em} \Xint-_{B_{i} }  |Dv|^{1+\frac{1}{s_{\varphi}}}  \omega(x) \, \mathrm{d}x \right)^{\frac{s_{\varphi}}{s_{\varphi}+1}}   \left(     \hspace{0.25em} \Xint-_{B_{i} }  \omega(x)^{-s_{\varphi}}  \, \mathrm{d}x  \right)^{\frac{1}{s_{\varphi}+1}}  \right) w(B_{i})    \\
    \notag &&  \leq   \sum_{i \in \mathbb{N} }  \varphi\left( \left( \frac{  [w]_{A_{1+ \frac{1}{s_{\varphi}}}}}{\omega(B_{i})} \int_{B_{i} }  |Dv|^{1+\frac{1}{s_{\varphi}}}  \omega(x) \, \mathrm{d}x \right)^{\frac{s_{\varphi}}{s_{\varphi}+1}}      \right) w(B_{i})    \\
    &&  \leq  c \sum_{i \in \mathbb{N} }  \int_{B_{i} }  \varphi\left( |Dv|\right)    \omega(x) \, \mathrm{d}x  \leq c \int_{E_{\lambda} }  \varphi\left( |Dv|\right)    \omega(x) \, \mathrm{d}x
 \end{eqnarray}
for some constant $c$ depending only on $n, s_{\varphi}$ and $
[w]_{A_{1+ \frac{1}{s_{\varphi}}}}$, where we have used the Jensen
inequality with a convex function
$\varphi(t^{\frac{s_{\varphi}}{s_{\varphi}+1}})$ for the first
inequality on the last line of \eqref{lipesti3}. Combining
\eqref{lipesti2} and \eqref{lipesti3}, we obtain the desired
estimate \eqref{lipsestidif}.
\end{proof}
\begin{remark}
    Applying  the Poincar\'e inequality together with the above Lemma \ref{liptrun} to a ball $B_{2\mathrm{diam}(\Omega)}$ with $\Omega \subset B_{2\mathrm{diam}(\Omega)}$, we obtain
    \begin{eqnarray}\label{rmkesti1}
        \notag  \hspace{-5mm}  \int_{\Omega} \varphi \left(  \frac{|v_\lambda|}{\mathrm{diam} (\Omega)}  \right)  \omega(x) \, \mathrm{d}x     \underset{ \eqref{poin1} }{\leq}  c \int_{\Omega} \varphi ( |  D v_\lambda| )  \omega(x) \, \mathrm{d}x
    \end{eqnarray}
for some constant $c$ depending on $n, s_{\varphi}, M$ and $[\omega]_{A_{1+\frac{1}{s_{\varphi}}}}$. Then it directly follows that $ \| v_\lambda  \|_{W_{\omega}^{1,\varphi}(\Omega)} \leq  c \| v  \|_{W_{\omega}^{1,\varphi}(\Omega)} $ for some constant  $c$ depending on $n, s_{\varphi}, M,  [\omega]_{A_{1+\frac{1}{s_{\varphi}}}} $ and $\mathrm{diam}(\Omega)$.
\end{remark}

We end this section with the following variant of Calder\'on-Zygmund-Krylov-Safonov decomposition lemma, which will be used later for the proof of weighted gradient estimates \eqref{mainesti3}.

\begin{lemma}\cite{MP}\label{calzyg}
    Assume that $\Omega \subset \mr^n$  is a $(M, r_{0})$-Lipschitz domain. and let  $\omega $ be an $ A_{s}$ weight for some $s>1$. Suppose that the sequence of balls $\{ B_{\rho_{0}} (y_{i})\}_{i=1}^{L}$ with centers $y_{i} \in \overline{\Omega}$ and a common radius $\rho_{0} \leq \frac{r_{0}}{4}$ covers $\Omega$. Let $C \subset D \subset \Omega$ be measurable sets for which there exists $0 < \epsilon < 1$ such that\newline
    \hspace*{2.7mm} \textup{(1)} $\omega(C) < \epsilon \omega(B_{\rho_{0} }(y_{i}))$ for all $i=1, 2, \cdots, L$, and \newline
    \hspace*{4mm}\textup{(2)}  for all $x \in \Omega$ and $\rho \in (0, 2\rho_{0} ]$, if $\omega(C \cap B_{\rho}(x)) \geq \epsilon \omega(B_{\rho}(x))$, then $B_{\rho}(x) \cap  \Omega \subset D$.
    Then we have the estimate
    $$
    \omega(C) \leq \epsilon c \omega(D).
    $$
    for some constant $c$ which depends only on $n, s, [\omega]_{A_{s}}$ and $ M$.
\end{lemma}

\section{Global weighted gradient estimates of very weak solutions}\label{Sec4}

From now on,  the letter $c$ represents a generic constant depending
on $n,s_{\varphi}, \nu, L$ and $M$ unless stated otherwise. To establish global weighted
gradient estimates, we employ a classical approach utilizing
Calder\'on-Zygmund decomposition as in \cite{MP}.  For the
comparison estimates which is the essence of the proof, we revisit
the earlier paper \cite{BL, BL2} to avoid duplication of
calculations. Here we only present the comparison estimates near the
boundary and refer to \cite[Lemma 5.1]{BL} for the interior
comparison estimates.

Consider a very weak solution $w$ to the following homogeneous
problem
\begin{equation}\label{refeq}
    \begin{cases}
        \mathrm{div\,}\mathcal{A}(x,Dw)=0 & \textrm{in}\ \Omega_{2r} \\
        w\in u+ W^{1,\varphi^{1-\sigma}}_{0}( \Omega_{2r} )
    \end{cases}
\end{equation}
for each $\Omega_{2r}=\Omega_{2r}(y_{0})$ with $y_{0} \in \partial \Omega $ and $r>0$, where $\sigma$ will be chosen below in Lemma \ref{compari}. Then we have the following unweighted comparison estimates.

\begin{lemma}\label{compari} \cite[Lemma 5.4]{BL2}
    Assume \eqref{youngcondi} and  \eqref{moncondi}. Suppose that $\Omega $ is a $(M, r_{0})$-Lipschitz domain. Then for any
    $\varepsilon \in (0,1]$, there exists a positive constant $\sigma = \sigma(n,
    s_{\varphi}, \nu, L, M) \leq \frac{1}{2} $ and $\delta = \delta(n,
    s_{\varphi}, \nu, L, M, \varepsilon ) $
    such that the following holds: Suppose that $u \in
    W_{0}^{1,\varphi^{1-\sigma}} (\Omega)$ is a very weak solution to \eqref{maineq3} with
    $\mathbf{f} \in W^{1,\varphi^{1-\sigma}}(\Omega)$. We further assume that  for $\Omega_{2r}=\Omega_{2r}(y_{0})$ with $y_{0} \in \partial \Omega $ and $r \in (0,r_{0}]$
    \begin{equation}\label{bdylambda}
      \notag   \Xint-_{B_{2r}}  [\varphi(|Du|)]^{1-\sigma} \, \mathrm{d}x \leq \Lambda,  \quad  \hspace{0.25em} \Xint-_{B_{2r}}  [ \varphi(|\mathbf{f}|)]^{1-\sigma} \, \mathrm{d}x \leq \delta\Lambda
    \end{equation}
    for some $\Lambda>0$. Then there exists a very weak solution $w
    \in W^{1,\varphi^{1-\sigma}}  ( \Omega_{2r}  ) $ to the equation
    \eqref{refeq} such that
    \begin{equation}\label{compa}
     \notag   \Xint-_{B_{2r}}|V(Du)-V(Dw)|^{2-2\sigma} \, \mathrm{d}x \leq \varepsilon \Lambda
    \end{equation}
    with the estimate
    \begin{equation}\label{compenergy}
     \notag   \Xint-_{B_{2r}} [\varphi(|Dw|)]^{1-\sigma} \, \mathrm{d}x  \leq  c \Lambda,  \quad   \left(  \hspace{0.25em} \Xint-_{B_{r}} [\varphi(|Dw|)]^{1+\sigma} \, \mathrm{d}x  \right)^{\frac{1-\sigma}{1+\sigma}}  \leq  c \Lambda,
    \end{equation}
    where the positive constants  $c$ and  $\sigma \in (0, \frac{1}{2}]$ depend  only on  $ n, s_{\varphi}, \nu, L$ and $M$.
\end{lemma}

To apply the above lemma to our setting, we change it into the following
 form of weighted estimates, where $\sigma \leq \frac{1}{2}$ is a universal
constant given in the lemmas just mentioned above.

\begin{proposition}\label{calzyg2}
    Assume \eqref{youngcondi} and  \eqref{moncondi}. Suppose that $\Omega $ is a $(M, r_{0})$-Lipschitz domain. Then for any
    $\varepsilon \in (0,1]$ and a weight $\omega \in RH_{1+\frac{2}{\sigma}} $, there exists a positive constant $\delta = \delta(n,
    s_{\varphi}, \nu, L, M, \varepsilon )   $
    such that the following holds: Suppose that $u  \in W_{0}^{1,\varphi^{1-\sigma}} (\Omega)$ is a very weak solution to \eqref{maineq3} with $\mathbf{f} \in W^{1,\varphi^{1-\sigma}}(\Omega)$. If there exists a number $\Lambda>0$ such that
    \begin{align}\label{calzyg2assumption}
        && \hspace{-5mm} \{\mathcal{M} ( \chi_{\Omega} \varphi (|D u|)^{1-\sigma}  ) (x) \leq \Lambda \}  \cap  \{   \mathcal{M} ( \chi_{\Omega} \varphi (\mathbf{|f|})^{1-\sigma}  ) (x) \leq \delta \Lambda \} \cap \Omega_{\rho}(y)  \neq \emptyset
    \end{align}
     for some ball with $\Omega_{\rho}(y)$ with $y \in \overline{\Omega}$ and $\rho \leq \frac{r_{0}}{8}$, then for any $T\geq 1$ we have the following decay estimates
    \begin{equation}\label{calzyg2result}
     \notag \frac{\omega(\{   \mathcal{M} ( \chi_{\Omega} \varphi (|D u|)^{1-\sigma}  ) (x) > T\Lambda \} \cap \Omega_{\rho}(y) )}{\omega(B_{\rho}(y))}  <  c_{*} (  {T}^{-\frac{1+\sigma/3}{1-\sigma}}  + \varepsilon{T}^{-\frac{2}{\sigma+2}}) ,
    \end{equation}
     where the constant $c_{*}$ depends only on $ n, s_{\varphi}, \nu, L, M$ and $[\omega]_{RH_{1+\frac{2}{\sigma}}}$.
\end{proposition}

\begin{proof}
    From \eqref{calzyg2assumption}, we see that there exists $x_{0} \in \Omega_{\rho}(y)  $ such that
    \begin{equation}\label{calzyg2pf1}
        \frac{1}{|B_{r}|}\int_{\Omega_{r}(x_{0})}   [\varphi(|Du|)]^{1-\sigma} \, \mathrm{d}x \leq \Lambda,  \quad  \frac{1}{|B_{r}|} \int_{\Omega_{r}(x_{0})} [ \varphi(|\mathbf{f}|)]^{1-\sigma} \, \mathrm{d}x \leq \delta\Lambda
    \end{equation}
     for any $r>0$. In particular, for any point $x \in \Omega_{\rho}(y)  $  we have
    \begin{equation}\label{calzyg2pf12}
         \mathcal{M} ( \chi_{\Omega}  \varphi (|D u|)^{1-\sigma}) (x) \leq \max \{ \mathcal{M} ( \chi_{\Omega_{2\rho}(y)} \varphi (|D u|)^{1-\sigma}) (x), 3^n\Lambda \}.
    \end{equation}
        Now we first consider the boundary case that $B_{4\rho}(y) \cap \partial\Omega \neq \emptyset$. Let $y_{0} \in \partial\Omega$ be a boundary point such that $|y-y_{0}|=dist(y, \partial\Omega)$ and let $w \in u +  W^{1,\varphi^{1-\sigma}}_{0}( \Omega_{8\rho}(y_{0}))$ be a very weak solution to the Dirichlet problem
        \begin{equation}\label{bdyrefeq}
           \notag \begin{cases}
                \mathrm{div\,}\mathcal{A}(x,Dw)=0 & \textrm{in}\  \Omega_{8\rho}(y_{0}), \\
                w= u & \textrm{on} \  \Omega_{8\rho}(y_{0}).
            \end{cases}
        \end{equation}
        Here we extend $u$ by zero to $\mr^n \backslash \Omega$ and then extend $w$ by $u$ to $\mr^n \backslash \Omega_{8\rho}(y_{0})$. Since
        $$
        B_{4\rho}(y) \subset B_{8\rho}(y_{0}) \subset B_{12\rho}(y) \subset B_{13\rho}(x_{0}),
        $$
        we immediately have that \eqref{calzyg2pf1} implies
    \begin{eqnarray} \label{compahypo}
        \notag&& \hspace{-5mm} \frac{1}{|B_{8\rho}|}\int_{\Omega_{8\rho}(y_{0})}   [\varphi(|Du|)]^{1-\sigma} \, \mathrm{d}x  \leq \left(\frac{13}{8} \right)^{n} \frac{1}{|B_{13\rho}|} \int_{\Omega_{13\rho}(x_{0})}  [\varphi(|Du|)]^{1-\sigma} \, \mathrm{d}x  \leq   \frac{13^{n}  \Lambda}{8^{n}}  \\
        \notag && \hspace{-5mm} \frac{1}{|B_{8\rho}|} \int_{\Omega_{8\rho}(y_{0})}  [ \varphi(|\mathbf{f}|)]^{1-\sigma} \, \mathrm{d}x \leq \left(\frac{13}{8} \right)^{n} \frac{1}{|B_{13\rho}|} \int_{\Omega_{13\rho}(x_{0})}   [ \varphi(|\mathbf{f}|)]^{1-\sigma} \, \mathrm{d}x \leq  \frac{13^{n}\delta\Lambda}{8^{n}}.
    \end{eqnarray}
  Then we are under the hypothesis of Lemma \ref{compari}. For any positive number $T \geq  3^n$, using \eqref{calzyg2pf12} and weak type (1,1) estimates for the maximal function we have
        \begin{eqnarray}\label{calzyg2pf2}
        \notag&& |\{ x \in \Omega_{\rho}(y) : \mathcal{M} ( \chi_{\Omega} \varphi (|D u|)^{1-\sigma}  )  > T\Lambda \}| \\
        \notag&& \hspace{5mm} \leq |\{ x \in \Omega_{\rho}(y) : \mathcal{M} ( \chi_{\Omega_{2\rho}(y)} \varphi (|D u|)^{1-\sigma}  )  > T\Lambda \}| \\
        \notag&& \hspace{5mm}  \leq \left| \left\{ x \in \Omega_{\rho}(y)  : \mathcal{M} ( \chi_{\Omega_{2\rho}(y)} \varphi (|D w|)^{1-\sigma}  )  > \frac{T\Lambda}{2c} \right\} \right| \\
        \notag&& \hspace{10mm} + \left| \left\{ x \in \Omega_{\rho}(y)  : \mathcal{M} ( \chi_{\Omega_{2\rho}(y)}|V(Du) - V(Dw)|^{2-2\sigma} ) > \frac{T\Lambda}{2c} \right\} \right| \\
        \notag&& \hspace{5mm} \leq \left(  \frac{c}{T\Lambda}\right)^{\frac{1+\sigma}{1-\sigma}} \int_{\Omega_{2\rho}(y)}   [\varphi(|Dw|)]^{1+\sigma} \, \mathrm{d}x \\
        \notag&& \hspace{15mm} + \frac{c}{T\Lambda} \int_{\Omega_{2\rho}(y)}   |V(Du) - V(Dw)|^{2-2\sigma} \, \mathrm{d}x \\
        &&  \hspace{5mm} \leq c_{1} (  {T}^{-\frac{1+\sigma}{1-\sigma}}  + \varepsilon{T}^{-1}) |B_{\rho}(y)|
    \end{eqnarray}
for some constant $c_{1}$ depending only on $ n, s_{\varphi}, \nu, L$ and $M$, where we have used Lemma \ref{compari} for the last inequality of \eqref{calzyg2pf2}.
For the interior case that  $B_{4\rho}(y) \subset \Omega$, let $w \in u +  W^{1,\varphi^{1-\sigma}}_{0}( B_{4\rho}(y))$ be a solution to the Dirichlet problem
\begin{equation}\label{intrefeq}
    \notag \begin{cases}
        \mathrm{div\,}\mathcal{A}(x,Dw)=0 & \textrm{in}\  B_{4\rho}(y), \\
        w= u & \textrm{on} \  B_{4\rho}(y).
    \end{cases}
\end{equation}
 Similarly as before, direct computation gives
\begin{eqnarray} \label{compahypo2}
    \notag&& \hspace{-5mm}
     \Xint-_{B_{4\rho}(y)}   [\varphi(|Du|)]^{1-\sigma} \, \mathrm{d}x  \leq \left(\frac{5}{4} \right)^{n}    \Xint-_{B_{5\rho}(x_{0})}  [\varphi(|Du|)]^{1-\sigma} \, \mathrm{d}x  \leq   \frac{5^{n}  \Lambda}{4^{n}}  \\
    \notag && \hspace{-5mm}  \Xint-_{B_{4\rho}(y)}  [ \varphi(|\mathbf{f}|)]^{1-\sigma} \, \mathrm{d}x \leq \left(\frac{5}{4} \right)^{n}   \Xint-_{B_{5\rho}(x_{0})}   [ \varphi(|\mathbf{f}|)]^{1-\sigma} \, \mathrm{d}x \leq  \frac{5^{n}\delta\Lambda}{4^{n}}
\end{eqnarray}
    and then we are under the hypothesis of Lemma 5.1 in \cite{BL}, the interior comparison estimates corresponding to Lemma \ref{compari}. Then the same type of inequality  \eqref{calzyg2pf2} follows for interior case, too. Therefore, for both interior and boundary case, we have the inequality
\begin{equation}\label{calzyg2pf4}
 |\{ x \in \Omega_{\rho}(y) : \mathcal{M} ( \chi_{\Omega} \varphi (|D u|)^{1-\sigma}  )  > T\Lambda \}| \leq c_{2} (  {T}^{-\frac{1+\sigma}{1-\sigma}}  + \varepsilon{T}^{-1}) |B_{\rho}(y)|,
\end{equation}
  where the constant $c_{2}$ depending only on $ n, s_{\varphi}, \nu, L$ and $M$. Note that if $T \geq 2 c_{2}$ the right hand side of \eqref{calzyg2pf4} is bounded above by $|B_{\rho}(y)|$. Since $\omega \in RH_{1+\frac{2}{\sigma}}$, we apply \eqref{aprhproperty} with $q={1+\frac{2}{\sigma}}$ to obtain
  \begin{eqnarray}\label{calzyg2pf5}
    \notag && \hspace{-5mm} \omega(\{ x \in \Omega_{\rho}(y) : \mathcal{M} ( \chi_{\Omega} \varphi (|D u|)^{1-\sigma}  )  > T\Lambda \}) \\
    \notag && \leq [\omega]_{RH_{1+\frac{2}{\sigma}}} \left[c_{2}  (  {T}^{-\frac{1+\sigma}{1-\sigma}}  + \varepsilon{T}^{-1})\right]^{\frac{2}{\sigma+2}} \omega(B_{\rho}(y)) \\
    \notag &&  <  c_{*} (  {T}^{-\frac{1+\sigma/3}{1-\sigma}}  + \varepsilon{T}^{-\frac{2}{\sigma+2}})   \omega(B_{\rho}(y))
  \end{eqnarray}
  for $T \geq \max \{ 3^n,  2c_{2} \}$, where the constant $c_{*}$ depending only on $ n, s_{\varphi}, \nu, L, M$ and $[w]_{RH_{1+\frac{2}{\sigma}}}$. For  $ 1 \leq T  \leq \max \{ 3^n, 2c_{2} \} $,
  if we replace $c_{*}$ with some constant greater than $\max \{ 3^n,  2c_{2}\}$, it follows that
  \begin{eqnarray}\label{calzyg2pf6}
    \notag && \hspace{-5mm} \omega(\{ x \in \Omega_{\rho}(y) : \mathcal{M} ( \chi_{\Omega} \varphi (|D u|)^{1-\sigma}  )  > T\Lambda \}) \leq   \omega(B_{\rho}(y)) \\
    \notag &&\hspace{20mm} <  c_{*} (  {T}^{-\frac{1+\sigma/3}{1-\sigma}}  + \varepsilon{T}^{-\frac{2}{\sigma+2}})   \omega(B_{\rho}(y))
  \end{eqnarray}
  since $c_{*}  {T}^{-\frac{1+\sigma/3}{1-\sigma}} > 1$. This completes the proof.
\end{proof}

We shortly remark that since the constant $T$ will be fixed as a large universal number
in the middle of the proof of Theorem \ref{mainthm2}, the lower bound for  $T$ in Proposition \ref{calzyg2} is not essential for our main proof. Combining Lemma \ref{calzyg} and Proposition \ref{calzyg2}, we get the following:
\begin{proposition}\label{calzyg3}
     Under the same assumption as in  Proposition \ref{calzyg2}, we further assume that $\omega \in A_{\frac{2}{2-\sigma}}$ and there exists a positive number $ \Lambda$  such that
    \begin{eqnarray}\label{calzyg3assumption}
     \notag && \omega( \{ x \in \Omega : \mathcal{M} ( \chi_{\Omega} \varphi (|D u|)^{1-\sigma}  ) (x) > T \Lambda \} ) \\
     && \hspace{40mm}  \leq  c_{*} (  {T}^{-\frac{1+\sigma/3}{1-\sigma}}  + \varepsilon{T}^{-\frac{2}{\sigma+2}}) \omega(B_{\frac{r_{0}}{8}}(y))
    \end{eqnarray}
    for any $y \in \overline{\Omega}$ and $T \geq 1$, where $c_{*}$ is the number given in Proposition \ref{calzyg2}. Then we have the following decay estimate
    \begin{eqnarray}\label{calzyg3result}
        \notag &&\hspace{-5mm}  \omega( \{ x \in \Omega : \mathcal{M} ( \chi_{\Omega} \varphi (|D u|)^{1-\sigma}  ) (x) > T \Lambda   \} )   \\
        \notag &&\hspace{5mm} \leq  c^{*} ({T}^{-\frac{1+\sigma/3}{1-\sigma}}  + \varepsilon{T}^{-\frac{2}{\sigma+2}})   \omega( \{ x \in \Omega : \mathcal{M} ( \chi_{\Omega} \varphi (|D u|)^{1-\sigma}  ) (x) > \Lambda  \} )      \\
        &&\hspace{20mm}   + c^{*} \omega( \{ x \in \Omega : \mathcal{M} ( \chi_{\Omega} \varphi (|\mathbf{f}|)^{1-\sigma}  ) (x) > \delta\Lambda  \} ) ,
    \end{eqnarray}
     for some constant $c^{*}$ depending only on $ n, s_{\varphi}, \nu, L, M, [w]_{A_{\frac{2}{2-\sigma}}}$ and $[w]_{RH_{1+\frac{2}{\sigma}}}$.
\end{proposition}

\begin{proof}
Set
$$
C :=  \{ x \in  \Omega : \mathcal{M} ( \chi_{\Omega} \varphi (|D u|)^{1-\sigma}  ) (x) > T\Lambda  \}
$$
and
    \begin{eqnarray}\label{Ddefinition}
        \notag &&\hspace{-5mm}  D :=  \{ x \in \Omega : \mathcal{M} ( \chi_{\Omega} \varphi (|D u|)^{1-\sigma}  ) (x) > \Lambda  \} \\
        \notag &&\hspace{5mm} \cup  \{ x \in  \Omega : \mathcal{M} ( \chi_{\Omega} \varphi (|D u|)^{1-\sigma}  ) (x) > \delta\Lambda  \} .
    \end{eqnarray}
For the case that $T > 4c_{*}^{2}$ so that $c_{*} (
{T}^{-\frac{1+\sigma/3}{1-\sigma}}  +
\varepsilon{T}^{-\frac{2}{\sigma+2}})< 1$, the assumption
\eqref{calzyg3assumption} gives the first hypothesis of Lemma
\ref{calzyg} with $\epsilon = c_{*} (
{T}^{-\frac{1+\sigma/3}{1-\sigma}}  +
\varepsilon{T}^{-\frac{2}{\sigma+2}})$. Moreover, suppose that there are $x \in \Omega$ and $\rho \in
(0, \frac{r_{0}}{8}]$ such that
$$
\omega(C \cap B_{\rho}(x)) \geq c_{*} (  {T}^{-\frac{1+\sigma/3}{1-\sigma}}  + \varepsilon{T}^{-\frac{2}{\sigma+2}}) \omega(B_{\rho}(x)).
$$
Then we apply Proposition \ref{calzyg2}, to discover that for any $z
\in B_{\rho}(x) \cap \Omega$ we have
\begin{eqnarray}\label{notin}
    \notag &&\hspace{-5mm}  z \notin   \Omega \cap \{ \mathcal{M} ( \chi_{\Omega} \varphi (|D u|)^{1-\sigma}  ) (x) \leq  \Lambda  \} \cap  \{ \mathcal{M} ( \chi_{\Omega} \varphi (|D u|)^{1-\sigma}  ) (x) \leq  \delta\Lambda  \} ,
\end{eqnarray}
which implies that $B_{\rho}(x) \cap \Omega \subset D$. This
verifies the second hypothesis of Lemma \ref{calzyg} and the desired
estimates \eqref{calzyg3result} follows. We replace $c^{*}$ with a
greater constant $\max  \{c^{*} , c_{*}^{3} \}  $ for the case that
$1 \leq T \leq 4c_{*}^{2} $. This completes the proof.
\end{proof}

We are now in a position to prove Theorem \ref{mainthm2}.
\begin{proof}[Proof of Theorem \ref{mainthm2}]

Let $\delta_{*} = \frac{\sigma}{2}$ and recall that $\mathbf{f} \in
L_{\omega}^{\varphi}(\Omega;\mr^n) \subset
L^{\varphi^{1-\delta_{*}}}(\Omega;\mr^n) $ for any weight $\omega
\in A_{\frac{1}{1-\delta_{*}}}$, where $\sigma$ is the positive
number given in Lemma \ref{compari}. For $\Lambda > 0 $ and for any
$y \in \overline{\Omega}$, weak $(1,1)$ estimates implies that
\begin{eqnarray}\label{weak11esti}
    \notag&& \hspace{-3mm}  |\{ x \in \Omega_{\frac{r_{0}}{8}}(y) : \mathcal{M} ( \chi_{\Omega} \varphi (|D u|)^{1-\sigma}  ) (x) > \Lambda \}  | \\
    \notag&&  \leq      |\{ x \in \Omega : \mathcal{M} ( \chi_{\Omega} \varphi (|D u|)^{1-\sigma}  ) (x) > \Lambda \}  |  \\
        \notag&&   \leq \frac{c}{\Lambda}\int_{\Omega}   [\varphi(|Du|)]^{1-\sigma}
        \, \mathrm{d}x    \underset{\eqref{msrden}}{\leq} \frac{c_{0} |\Omega_{\frac{r_{0}}{8}}(y)|}{\Lambda r_{0}^{n}}\int_{\Omega}   [\varphi(|Du|)]^{1-\sigma}
        \, \mathrm{d}x
    \end{eqnarray}
  for some constant $c_{0}$ depending only on $ n $ and $ M$. Choose
  \begin{eqnarray}\label{lambdazero}
    \Lambda_{0} :=  \frac{4c_{0} }{r_{0}^n} \int_{\Omega} \left\{ [\varphi(|Du|)]^{1-\sigma}
    + \frac{ [ \varphi(|\mathbf{f}|)]^{1-\sigma}}{\delta} \right\} \, \mathrm{d}x,
  \end{eqnarray}
   where $\delta$ is given in Proposition \ref{calzyg2}. Then for any  $\Lambda > \Lambda_{0} $ we have
   \begin{eqnarray}\label{dumeasure}
    \notag    |\{ x \in \Omega_{\frac{r_{0}}{8}}(y) : \mathcal{M} ( \chi_{\Omega} \varphi (|D u|)^{1-\sigma}  ) (x)>\Lambda \}  |  \leq \frac{ |\Omega_{\frac{r_{0}}{8}}(y)|}{4}
   \end{eqnarray}
    and similarly we have
    \begin{eqnarray}\label{fmeasure}
     \notag   |\{ x \in \Omega_{\frac{r_{0}}{8}}(y) : \mathcal{M} ( \chi_{\Omega} \varphi (|D u|)^{1-\sigma}  ) (x) > \delta \Lambda \}  |  \leq \frac{ |\Omega_{\frac{r_{0}}{8}}(y)|}{4},
    \end{eqnarray}
   which leads to the assumption \eqref{calzyg2assumption} with $\rho = \frac{r_{0}}{8}$. Applying Proposition \ref{calzyg2}, we obtain the estimate \eqref{calzyg3assumption} and we are under the hypothesis of Proposition \ref{calzyg3}. Define the upper level sets as
$$
\mathcal{E}_{\Lambda} := \{x \in \Omega :  \mathcal{M} (
\chi_{\Omega} \varphi (|D u|)^{1-\sigma})  >\Lambda \}  \ \ \
$$
and
$$
\mathcal{E}^{f}_{\Lambda} := \{x \in  \Omega   :  \mathcal{M} ( \chi_{\Omega} \varphi (|\mathbf{f}|)^{1-\sigma}) >\Lambda \}
$$
for $\Lambda> 0$. Then the inequality $\eqref{calzyg3result}$ becomes
\begin{eqnarray}\label{elambdaestimate}
    \notag &&\hspace{-5mm}  \omega( \mathcal{E}_{T\Lambda} )    \leq  c(n) ({T}^{-\frac{1+\sigma/3}{1-\sigma}}  + \varepsilon{T}^{-\frac{2}{\sigma+2}})   \omega( \mathcal{E}_{\Lambda} )  + c(n) \omega( \mathcal{E}^{f}_{\delta\Lambda} ).
\end{eqnarray}
We next introduce truncation functions
$$
\left[  \mathcal{M} ( \chi_{\Omega} \varphi (|D u|)^{1-\sigma}) \right]_{t}:=\min \{ \mathcal{M} ( \chi_{\Omega} \varphi (|D u|)^{1-\sigma}) , t \} \quad (t>0) .
$$
For $T>1$, we define $ t_{0} := T\Lambda_0$. Then for $t > t_{0}$ Fubini's theorem gives
\begin{eqnarray}\label{fubiniinteg1}
\nonumber &&   \hspace{-5mm} \int_{\Omega}  \left[  \mathcal{M} ( \chi_{\Omega} \varphi (|D u|)^{1-\sigma}) \right]_{t}^{\frac{1}{1-\sigma}} \omega(x) \, \mathrm{d}x  \\
\nonumber &&=  \frac{1}{1-\sigma}   \int_0^{t} \Lambda^{\frac{\sigma}{1-\sigma}} \omega (\mathcal{E}_{\Lambda})   \, d\Lambda \\
\nonumber && \leq \frac{1}{1-\sigma} \int_0^{t_{0}} \Lambda^{\frac{\sigma}{1-\sigma}}\omega(\mathcal{E}_{\Lambda})   \, d\Lambda   +   \frac{1}{1-\sigma}  \int_{t_{0}}^{t} \Lambda^{\frac{\sigma}{1-\sigma}} \omega(\mathcal{E}_{\Lambda})  \, d\Lambda  \\
\nonumber &&  \hspace{-2.5mm}  \underset{\eqref{lambdazero}}{\leq}  T^{\frac{1}{1-\sigma}} \Lambda_{0}^{\frac{1}{1-\sigma}}  \omega(\Omega) \\
\nonumber && \hspace{5mm} + \frac{c^{*} }{1-\sigma} \left( {T}^{-\frac{1+\sigma/3}{1-\sigma}}  + \varepsilon{T}^{-\frac{2}{\sigma+2}}  \right)       \int_{t_{0}}^{t} \Lambda^{\frac{\sigma}{1-\sigma}} \omega(\mathcal{E}_{\Lambda/T})  \, d\Lambda   \\
&& \hspace{25mm}   +   \frac{c^{*} }{1-\sigma}  \int_{t_0}^{t}  \Lambda^{\frac{\sigma}{1-\sigma}} \omega(\mathcal{E}^{f}_{\delta\Lambda/T}) \, d\Lambda.
\end{eqnarray}
Using change of variable, we can estimate each integral appeared in the last line as follows:
\begin{eqnarray}\label{fubiniinteg2}
    \nonumber &&   \hspace{-15mm}  \int_{t_{0}}^{t} \Lambda^{\frac{\sigma}{1-\sigma}} \omega(\mathcal{E}_{\Lambda/T})  \, d\Lambda = T^{\frac{1}{1-\sigma}}
    \int_{t_{0}/T}^{t/T} \Lambda^{\frac{\sigma}{1-\sigma}} \omega(\mathcal{E}_{\Lambda})  \, d\Lambda \\
    \nonumber && \hspace{-10mm} \leq T^{\frac{1}{1-\sigma}}
    \int_{0}^{t} \Lambda^{\frac{\sigma}{1-\sigma}} \omega(\mathcal{E}_{\Lambda})  \, d\Lambda  \\
    && \hspace{-10mm} \leq \left(1-\sigma\right)  T^{\frac{1}{1-\sigma}}
       \int_{\Omega}  \left[  \mathcal{M} ( \chi_{\Omega} \varphi (|D u|)^{1-\sigma}) \right]_{t}^{\frac{1}{1-\sigma}} \omega(x) \, \mathrm{d}x
\end{eqnarray}
and
\begin{eqnarray}\label{fubiniinteg3}
    \nonumber &&   \hspace{-15mm} \int_{t_0}^{t}  \Lambda^{\frac{\sigma}{1-\sigma}} \omega(\mathcal{E}^{f}_{\delta\Lambda/T}) \, d\Lambda \\
    && \hspace{-10mm} \leq \frac{\left(1-\sigma\right)  T^{\frac{1}{1-\sigma}}}{\delta^{\frac{1}{1-\sigma}}}   \int_{\Omega}  \left[  \mathcal{M} ( \chi_{\Omega} \varphi (|\mathbf{f}|)^{1-\sigma}) \right]^{\frac{1}{1-\sigma}} \omega(x) \, \mathrm{d}x .
\end{eqnarray}
Combining \eqref{fubiniinteg1}, \eqref{fubiniinteg2} and \eqref{fubiniinteg3}, we get
\begin{eqnarray}\label{fubiniinteg4}
    \nonumber && \hspace{-15mm}  \int_{\Omega}  \left[  \mathcal{M} ( \chi_{\Omega} \varphi (|D u|)^{1-\sigma}) \right]_{t}^{\frac{1}{1-\sigma}} \omega(x) \, \mathrm{d}x \\
    \nonumber && \hspace{-10mm} \leq  T^{\frac{1}{1-\sigma}} \Lambda_{0}^{\frac{1}{1-\sigma}}  \omega(\Omega)   \\
    \nonumber && \hspace{-5mm} + c^{*} \left( {T}^{-\frac{\sigma/3}{1-\sigma}}  + \varepsilon{T}^{\frac{1}{1-\sigma} -\frac{2}{\sigma+2}}  \right)  \int_{\Omega}  \left[  \mathcal{M} ( \chi_{\Omega} \varphi (|\mathbf{f}|)^{1-\sigma}) \right]_{t}^{\frac{1}{1-\sigma}} \omega(x) \, \mathrm{d}x  \\
    \nonumber &&   + \frac{ c^{*}  T^{\frac{1}{1-\sigma}}  }{\delta^{\frac{1}{1-\sigma}}}  \int_{\Omega}  \left[  \mathcal{M} ( \chi_{\Omega} \varphi (|\mathbf{f}|)^{1-\sigma}) \right]_{t}^{\frac{1}{1-\sigma}} \omega(x) \, \mathrm{d}x.
\end{eqnarray}
 Now we choose $T= \left( 4 c^{*}\right )^{\frac{3}{\sigma}}$ and $\varepsilon = \frac{1}{4c^{*} T^{2}}$ in order to have the following estimate:
\begin{equation}\label{fubiniinteg5}
   \nonumber \int_{\Omega}  \left[  \mathcal{M} ( \chi_{\Omega} \varphi (|D u|)^{1-\sigma}) \right]_{t}^{\frac{1}{1-\sigma}} \omega(x) \,  \mathrm{d}x     \leq c \Lambda_{0}^{\frac{1}{1-\sigma}}  \omega(\Omega)+ c  \int_{\Omega}\varphi(|\mathbf{f}|)\omega(x) \, \mathrm{d}x  .
\end{equation}
We remark that $\delta$ now depends only on  $ n, s_{\varphi},
\nu, L, M, [\omega]_{A_{\frac{1}{1-\delta_{*}}}}$ and $  [\omega]_{RH_{1+\frac{1}{\delta_{*}}}}$. Therefore, using the unweighted a priori estimates, we have
$$
\Lambda_{0} = \frac{4c_{0} }{r_{0}^n}  \int_{\Omega} \left\{ [\varphi(|Du|)]^{1-\sigma}
+ \frac{ [ \varphi(|\mathbf{f}|)]^{1-\sigma}}{\delta} \right\} \, \mathrm{d}x \leq \frac{c}{r_{0}^{n}} \int_{\Omega} [ \varphi(|\mathbf{f}|)]^{1-\sigma}  \, \mathrm{d}x
$$
for some constant $c$ depending only on $ n, s_{\varphi},
\nu, L, M, [\omega]_{A_{\frac{1}{1-\delta_{*}}}}$ and $  [\omega]_{RH_{1+\frac{1}{\delta_{*}}}}$.

Take a ball $B$ such that $B:=B_{\mathrm{diam}(\Omega)} \supseteq
\Omega$. By the H\"older inequality, we obtain
\begin{eqnarray}\label{mainproof1}
\nonumber &&   \hspace{-1mm}  \int_{\Omega}\left[  \mathcal{M} ( \chi_{\Omega} \varphi (|D u|)^{1-\sigma}) \right]_{t}^{\frac{1}{1-\sigma}}  \omega(x) \, \mathrm{d}x  \\
\nonumber &&  \hspace{-6mm}  \leq c \omega(\Omega) \left(  \frac{c}{r_{0}^n} \int_{\Omega}  [\varphi(|\mathbf{f}|)]^{1-\sigma}
 \, \mathrm{d}x \right)^{\frac{1}{1-\sigma}}   + c \int_{\Omega} \varphi(|\mathbf{f}|)   \omega(x) \, \mathrm{d}x \\
 \nonumber &&  \hspace{-6mm}  \leq c   \omega(\Omega) \left( \frac{\mathrm{diam}(\Omega) }{r_{0}}  \right)^{\frac{n}{1-\sigma}} \left(   \, \Xint-_{B}   \varphi(|\mathbf{f}|) \omega
 \, \mathrm{d}x \right) \left(   \, \Xint-_{B}   \omega^{\frac{\sigma-1}{\sigma}}
 \, \mathrm{d}x \right) ^{\frac{\sigma}{1-\sigma}}   + c \int_{\Omega} \varphi(|\mathbf{f}|)   \omega(x) \, \mathrm{d}x \\
 \nonumber &&  \hspace{-6mm}  \leq c  \omega(\Omega) [\omega]_{A_{\frac{1}{1-\sigma}}} \left( \frac{\mathrm{diam}(\Omega) }{r_{0}}  \right)^{\frac{n}{1-\sigma}}   \left(  \frac{1}{\omega(B)}  \int_{B}   \varphi(|\mathbf{f}|) \omega
 \, \mathrm{d}x \right)    + c \int_{\Omega} \varphi(|\mathbf{f}|)   \omega(x) \, \mathrm{d}x  \\
 \nonumber &&  \hspace{-6mm}  \leq   c \left[   \left( \frac{\mathrm{diam}(\Omega) }{r_{0}}  \right)^{\frac{n}{1-\sigma}}  + 1 \right]  \int_{\Omega} \varphi(|\mathbf{f}|)   \omega(x) \, \mathrm{d}x.
\end{eqnarray}
Letting $t \rightarrow \infty$, we obtain the desired estimate
\begin{eqnarray}\label{mainproof2}
 \nonumber \int_{\Omega}   \varphi(|Du|)  \omega(x) \, \mathrm{d}x   \leq c  \int_{\Omega}  \varphi(|\mathbf{f}|)  \omega(x) \, \mathrm{d}x
\end{eqnarray}
for some positive constant $c$ depending only on  $ n, s_{\varphi},
\nu, L, M, [\omega]_{A_{\frac{1}{1-\delta_{*}}}},
[\omega]_{RH_{1+\frac{1}{\delta_{*}}}}$ and
$\frac{\mathrm{diam}(\Omega)}{r_{0}}$. This completes the
proof.
\end{proof}

\section{Existence of very weak solution}\label{Sec5}
This section is devoted to proving Theorem 2.1, the second main result of our paper.
Following the idea of the proof for existence of very weak solution described
in \cite{BDS, BS}, we first start with the following weighted div-curl lemma.
\begin{theorem}\label{divcurl}
Assume \eqref{youngcondi} and suppose that  $\Omega$ is a bounded Lipschitz domain.
Then for a given weight $\omega \in \mathcal{A}_{1+\frac{1}{s_{\varphi}}}$, assume that a sequence of vector-valued measurable functions $(a_{k}, b_{k})_{k=1}^{\infty} : \Omega \rightarrow \mr^n \times \mr^n$ satisfy
\begin{eqnarray}\label{divcurlasm1}
    \notag      a_{k} \rightharpoonup  a  &   \textrm{weakly}  \,    \textrm{in} \, \, L_{\omega}^{\varphi}(\Omega)    \\
    \notag  b_{k} \rightharpoonup  b    &   \textrm{weakly}  \,    \textrm{in} \, \, L_{\omega}^{\varphi^{*}}(\Omega)
\end{eqnarray}
for some $a \in L_{\omega}^{\varphi}(\Omega)  $ and $ b \in L_{\omega}^{\varphi^{*}}(\Omega)   $. We further assume that
\begin{equation}\label{divcurlasm2}
      \lim_{k \rightarrow \infty} \int_{\Omega} b_{k} \cdot D  \eta_{k} \, \mathrm{d}x = 0
\end{equation}
 for any sequence $\eta_{k} \in C_{0}^{\infty}(\Omega)$ such that
 \begin{equation}\label{divcurlasm3}
    D \eta_{k} \rightharpoonup 0  \quad \textrm{weakly} \, \,    \textrm{in} \, \, (L_{\omega}^{\varphi^{*}}(\Omega))^{*} ,
 \end{equation}
and  there holds
\begin{equation}\label{divcurlasm4}
\lim_{k \rightarrow \infty} \int_{\Omega} a_{k}^{i} \partial_{j}\xi_{k} -  a_{k}^{j} \partial_{i}\xi_{k}   \, \mathrm{d}x = 0  \, \, \,   \textrm{for all} \, \, \, \,  i, j = 1, \cdots , n.
\end{equation}
for any sequence $\xi_{k} \in C_{0}^{\infty}(\Omega)$ such that
\begin{equation}\label{divcurlasm5}
\notag  D \xi_{k} \rightharpoonup 0  \quad \textrm{weakly} \, \,    \textrm{in} \, \, (L_{\omega}^{\varphi}(\Omega))^{*} ,
\end{equation}
Then we have
\begin{equation}\label{divcurlresult}
  \notag   a_{k} \cdot b_{k} \, \omega \rightarrow  a \cdot b \,  \omega
\end{equation}
in the sense of distributions in $\Omega$.
\end{theorem}

\begin{proof}
    Take a ball $B$ such that $ \Omega \subset\subset  B$ and extend  $a_{k} $ by zero outside $\Omega$. Since $a_{k} \in  L_{\omega}^{\varphi}(B) \subset L^{q}(B) $  for some $q>1$, we can find the unique solution $d^{k}$ to Laplace equation
    \begin{equation}\label{divproof1}
        \begin{cases}
            \Delta d_{k}=a_{k} \quad    & \textrm{in}\ B  \\
            u=0 & \textrm{on}\ \partial B.
        \end{cases}
    \end{equation}
    Then it directly follows that
    \begin{equation}\label{divproof11}
           d_{k} \rightharpoonup d \quad \textrm{weakly in} \, \,  W^{2,q} (B; \mr^n),
    \end{equation}
    for some $ d \in  W^{2,q} (B; \mr^n) $, which gives
    \begin{equation}\label{divproof12}
      \notag  \Delta d=a \quad  \textrm{in}\ B.
    \end{equation}
 On the other hand, the  Calder\'on-Zygmund theory for Laplace equation underlined in weighted Orlicz spaces  gives
    \begin{equation}\label{divproof2}
      \notag  \int_{B}   \varphi(|D^2d_{k}|) \omega \, \mathrm{d}x   \leq c  \int_{B}  \varphi(|a_{k}|) \omega \, \mathrm{d}x \leq
      C
    \end{equation}
for some constant $C$ independent of $k$, see \cite{BR}. This
implies that
\begin{equation}\label{divproof3}
  \notag  d_{k} \rightharpoonup d  \quad \textrm{weakly} \, \,    \textrm{in} \, \, W_{\omega}^{2,\varphi}(B) ,
\end{equation}
for some $d \in W_{\omega}^{2,\varphi}(B)$.
Now, it is enough to check that
\begin{equation}\label{divproof5}
    b_{k} \cdot (a_{k}-D\, \mathrm{div} d_{k}) \omega \rightarrow  b \cdot (a - D\, \mathrm{div} d) \omega
\end{equation}
and that
\begin{equation}\label{divproof52}
b_{k} \cdot (D\, \mathrm{div} d_{k}) \omega \rightarrow b \cdot ( D\, \mathrm{div} d) \omega
\end{equation}
in the sense of distributions.

By taking the curl of both sides to the equality \eqref{divproof1}, we observe that standard Calder\'on-Zygmund theory with the assumption \eqref{divcurlasm4} gives the convergence
\begin{equation}\label{divproof522}
    \notag \partial_{i}d_{k}^{j} - \partial_{j}d_{k}^{i}  \rightarrow \partial_{i}d_{k}^{j} - \partial_{j}d_{k}^{i} \, \, \quad \textrm{strongly in} \, \,  W_{\omega}^{1,\varphi} (\Omega; \mr^n).
\end{equation}
 Combining this convergence with the following equality
\begin{equation}\label{divproof53}
    \notag  a_{k}^{j}- \partial_{j} \mathrm{div} d_{k} = \sum_{m=1}^{n} \partial^{2}_{mm} d_{k}^{j} - \partial_{j}\partial_{m}d_{k}^{m} = \sum_{m=1}^{n} \partial_{m} \left( \partial_{m}d_{k}^{j} - \partial_{j}d_{k}^{m} \right),
\end{equation}
we can conclude that the sequence $ a_{k}-D\, \mathrm{div} d_{k} $ strongly converges to $ a -D\, \mathrm{div} d  $ in $L_{\omega}^{\varphi}(\Omega)$. Therefore, we obtain
\begin{eqnarray}\label{divproof6}
    \notag&&\hspace{-15mm} \lim_{k \rightarrow \infty} \int_{\Omega} b_{k} \cdot (a_{k}-D\, \mathrm{div} d_{k}) \omega \eta \, \mathrm{d}x  \\
    \notag&&\hspace{-5mm}= \lim_{k \rightarrow \infty } \int_{\Omega}  b_{k} \cdot (a -D\, \mathrm{div} d ) \omega \eta \, \mathrm{d}x  \\
    \notag&& + \lim_{k \rightarrow \infty } \int_{ \Omega}   b_{k} \cdot (a_{k} -a -D\, \mathrm{div} d_{k} +  D\, \mathrm{div} d ) \omega \eta \, \mathrm{d}x   \\
    \notag&&\hspace{-5mm}=  \int_{\Omega}  b  \cdot (a -D\, \mathrm{div} d ) \omega \eta \, \mathrm{d}x  \\
    \notag && + \lim_{k \rightarrow \infty } \int_{ \Omega}  b_{k} \cdot (a_{k} -a -D\, \mathrm{div} d_{k} +  D\, \mathrm{div} d ) \omega \eta \, \mathrm{d}x
\end{eqnarray}
for any $\eta \in C_{0}^{\infty}(\Omega)$, where the second term of the last equality can be estimated as follows:
\begin{eqnarray}\label{divproof7}
   \notag && \left|  \,  \int_{ \Omega}   b_{k} \cdot (a_{k} -a -D\, \mathrm{div} d_{k} +  D\, \mathrm{div} d ) \omega \eta \, \mathrm{d}x \, \right|  \\
   && \hspace{5mm}  \leq c || \eta ||_{\infty} ||  a_{k} -a -D\, \mathrm{div} d_{k} +  D\, \mathrm{div} d ||_{L^{\varphi}_{\omega}(\Omega)} || b_{k} ||_{L^{\varphi^{*}}_{\omega}(\Omega)}.
\end{eqnarray}
Since the strong convergence of $ a_{k}-D\, \mathrm{div} d_{k} $ implies that the last expression \eqref{divproof7} converges to $0$ as $k \rightarrow \infty$, the convergence \eqref{divproof5} follows.

To prove the convergence \eqref{divproof52}, we denote
$$
e_{k} := \upsilon \, \mathrm{div} \, d_{k} \,   \quad \textrm{and}
\quad e := \upsilon \,  \mathrm{div} \, d   \, ,
$$
where $ \upsilon    \in C_{0}^{\infty}(B) $ is a fixed function
which is identically one in $\Omega$. Then the condition
\eqref{divcurlasm3} and the convergence \eqref{divproof11} imply
that
\begin{eqnarray}\label{divproof71}
   \notag && e_{k} \rightharpoonup e  \quad \textrm{weakly} \, \,    \textrm{in} \, \, W_{0}^{1,q}(B),  \\
    &&  D e_{k} \rightharpoonup  D e  \quad \textrm{weakly} \, \,    \textrm{in} \, \, L_{\omega}^{\varphi}(B).
\end{eqnarray}
Applying Lemma \ref{liptrun} to $e_{k}$ with a fixed constant
$\lambda
> 0 $, we can find a sequence $e_{k,\lambda}$ such that
$$e_{k,\lambda}(x)=e_{k}(x),   \quad  De_{k,\lambda}(x)=De_{k}(x) \quad    \quad \mathrm{a.e.} \ x \in B \backslash E_{\lambda}  $$
with the estimates $\varphi(|De_{k,\lambda}|) \leq c \lambda $  for
a.e. $x \in B$. Then using Lemma \ref{liptrun}, after passing
to a subsequence, we can find a Lipschitz function $e_{\lambda}$
such that
\begin{eqnarray}\label{divproof72}
    \notag && De_{k,\lambda} \rightharpoonup De_{\lambda}  \quad \textrm{weakly}^{*} \, \,    \textrm{in} \, \, L^{\infty}(B),  \\
    \notag &&  D e_{k,\lambda} \rightharpoonup  D e_{\lambda}  \quad \textrm{weakly} \, \,    \textrm{in} \, \, L_{\omega}^{\varphi}(B), \\
    &&   e_{k,\lambda} \rightarrow  e_{\lambda}  \quad \textrm{strongly} \, \,    \textrm{in} \, \, C(B).
\end{eqnarray}
Combining \eqref{divcurlasm2} with \eqref{divproof72}, we have
\begin{eqnarray}\label{divproof73}
    \notag && \hspace{-5mm}  \lim_{k \rightarrow \infty} \int_{\Omega} b^{k} \cdot D e_{k,\lambda} \eta \, \mathrm{d}x  =  \lim_{k \rightarrow \infty} \int_{\Omega} b^{k} \cdot (D e_{k,\lambda} - D e_{\lambda}) \eta \, \mathrm{d}x  +  \int_{\Omega} b \cdot D e_{\lambda} \eta \, \mathrm{d}x  \\
    \notag && \hspace{-7mm} = \lim_{k \rightarrow \infty} \int_{\Omega} b^{k} \cdot  D[ ( e_{k,\lambda} -  e_{\lambda}) \eta ] \, \mathrm{d}x  - \lim_{k \rightarrow \infty} \int_{\Omega} b^{k} \cdot ( e_{k,\lambda} -  e_{\lambda}) D\eta \, \mathrm{d}x  +  \int_{\Omega} b  \cdot D e_{\lambda} \eta \, \mathrm{d}x \\
    \notag && \hspace{-7mm}  =   \int_{\Omega} b \cdot D e_{\lambda} \eta \, \mathrm{d}x
\end{eqnarray}
for any $\eta \in C_{0}^{\infty}(\Omega)$. Then it follows that
\begin{equation}\label{divproof74}
   \notag  b_{k} \cdot D e_{k,\lambda} \rightharpoonup   b \cdot D e_{\lambda}  \quad \textrm{weakly} \, \,    \textrm{in} \, \, L^{1}(\Omega).
\end{equation}
To show the convergence \eqref{divproof52}, we estimate
\begin{eqnarray}\label{divproof8}
    \notag&&\hspace{-7mm}  \lim_{k \rightarrow \infty}  \left|  \int_{\Omega} [ b_{k} \cdot (D\, \mathrm{div} d_{k}) - b \cdot (D\, \mathrm{div} d)] \omega \eta \, \mathrm{d}x  \right|  = \lim_{k \rightarrow \infty}  \left|  \int_{\Omega} (b_{k} \cdot D e_{k} - b \cdot D e )  \omega \eta \, \mathrm{d}x  \right|  \\
    \notag&&\hspace{-7mm} \leq   \lim_{k \rightarrow \infty}  \left|  \int_{\Omega} (b_{k} \cdot D e_{k,\lambda} - b \cdot D e_{k} )  \omega \eta \, \mathrm{d}x  \right| +  \| \eta \|_{\infty} \limsup_{k \rightarrow \infty}   \int_{\Omega} |b_{k}| | D( e_{k,\lambda}-  e_{k}) |  \omega  \, \mathrm{d}x   \\
    \notag&& \hspace{-2mm}  +     \left|  \int_{\Omega} b  \cdot D( e_{\lambda} - e )  \omega \eta \, \mathrm{d}x  \right|  \\
    \notag&&\hspace{-7mm} \leq  \lim_{k \rightarrow \infty}  \left|  \int_{\Omega} \frac{(b_{k} \cdot D e_{k,\lambda} - b \cdot D e_{\lambda} )  \omega \eta }{1+ \varepsilon \omega} \, \mathrm{d}x  \right| +  \| \eta \|_{\infty} \limsup_{k \rightarrow \infty}   \int_{\Omega} |b_{k}| | D( e_{k,\lambda}-  e_{k}) |  \omega  \, \mathrm{d}x \\
    \notag&&\hspace{-2mm} + \| \eta \|_{\infty} \limsup_{k \rightarrow \infty}   \int_{\Omega} \frac{ \varepsilon \omega^2   (|b_{k}| | D e_{k,\lambda}| + |b| |D e_{\lambda}| )   }{1+ \varepsilon \omega} \, \mathrm{d}x  + \left|  \int_{\Omega} b  \cdot D( e_{\lambda} - e )  \omega \eta \, \mathrm{d}x  \right|    \\
    \notag&&\hspace{-7mm} \leq  \| \eta \|_{\infty} \limsup_{k \rightarrow \infty}   \int_{\Omega} \frac{ \varepsilon \omega^2  |b_{k}| | D e_{k,\lambda}|    }{1+ \varepsilon \omega} \, \mathrm{d}x   + \| \eta \|_{\infty} \limsup_{k \rightarrow \infty}   \int_{\Omega} \frac{ \varepsilon \omega^2    |b| |D  e_{\lambda}| }{1+ \varepsilon \omega} \, \mathrm{d}x  \\
    \notag&&\hspace{-2mm} + \| \eta \|_{\infty} \limsup_{k \rightarrow \infty}   \int_{\Omega} |b_{k}| | D( e_{k,\lambda}-  e_{k}) |  \omega  \, \mathrm{d}x   + \left|  \int_{\Omega} b  \cdot D( e_{\lambda} - e )  \omega \eta \, \mathrm{d}x  \right| \\
    \notag && \hspace{-7mm} =: I_{1}+ I_{2}+ I_{3}+ I_{4}
\end{eqnarray}
for arbitrary $\varepsilon>0$, where we have used the fact that $\frac{\omega\eta}{1+ \varepsilon \omega }$ is a bounded function with \eqref{divproof71} for the last inequality. Then it remains to prove that each $I_{j}$($j=1,2,3,4$) vanishes when $\varepsilon \rightarrow 0$ and $\lambda \rightarrow \infty$.  Convergence of $\lim_{\varepsilon \rightarrow 0}  I_{2} =0 $ directly follows from the Lebesgue convergence theorem with Lemma \ref{liptrun}. To estimate $I_{1}$, we compute
\begin{eqnarray}\label{divproof81}
    \notag&&\hspace{-7mm}   \limsup_{ \varepsilon \rightarrow 0}  I_{1}  \leq    \| \eta \|_{\infty}  \limsup_{ \varepsilon \rightarrow 0} \limsup_{k \rightarrow \infty}  \left|  \int_{\Omega} \frac{ \varepsilon \omega^2   |b_{k}| | D e_{k, \lambda}|    }{1+ \varepsilon \omega} \, \mathrm{d}x  \right|  \\
    \notag &&\leq C \limsup_{\varepsilon \to 0}\limsup_{k\to \infty}(\varphi^{*})^{-1} \left(\int_{\Omega\cap\{w>\lambda\}}\varphi^{*}({|b_k|})\omega \frac{\varepsilon \omega}{1+\varepsilon \omega} \, \mathrm{d}x \right)\\
    \notag&&\quad+C\limsup_{\varepsilon \to 0}\limsup_{k\to \infty} (\varphi^{*})^{-1}\left(\int_{\Omega\cap\{w\le \lambda\}}\varphi^{*}({|b_k|})\omega \frac{\varepsilon \omega}{1+\varepsilon \omega} \, \mathrm{d}x \right) \\
    \notag &&  \leq C\limsup_{k\to \infty}(\varphi^{*})^{-1} \left(\int_{\Omega\cap\{w>\lambda\}}\varphi^{*}({|b_k|})\omega \, \mathrm{d}x \right).
\end{eqnarray}
Since $|\{\omega >\lambda\}|\le C/\lambda$, we let $\lambda \to \infty$ in the last inequality to deduce
\begin{equation}    \label{lim2}
      \notag   \limsup_{\lambda\to\infty}\limsup_{\varepsilon \to 0} \limsup_{k \rightarrow \infty}   \int_{\Omega} \frac{ \varepsilon \omega^2  |b_{k}| | D e_{k,\lambda}|    }{1+ \varepsilon \omega} \, \mathrm{d}x  =0.
\end{equation}
Then we estimate $I_{3}$. By the H\"{o}lder inequality, we have
\begin{equation}
    \notag  \begin{split}
        \limsup_{\lambda \to \infty}I_{3}&=C\limsup_{\lambda \to \infty}\limsup_{k\to \infty}\int_{\Omega}|b_k| |D (e_k-e_{k,\lambda})|\omega \, \mathrm{d}x \\
        &=C\limsup_{\lambda \to \infty}\limsup_{k\to \infty}\int_{\Omega\cap\{M(D e_k)>\lambda\}}|b_k| |D(e_k-e_{k,\lambda})|\omega \, \mathrm{d}x \\
        &\le C\limsup_{\lambda \to \infty}\limsup_{k\to \infty} (\varphi^{*})^{-1}\left(\int_{\Omega\cap\{M(D e_k)>\lambda\}} \varphi^{*}({|b_k|})\omega  \, \mathrm{d}x \right) =0,
    \end{split}
    \label{lim4}
\end{equation}
where the last inequality follows from the fact that $|\{M(D
e_k)>\lambda\}|\le C/\lambda$. Finally, it suffices to show
\begin{equation}
    \begin{split}
        \lim_{\lambda \to \infty}I_{4}=&\lim_{\lambda \to \infty}\left|\int_{\Omega}b \cdot D (e-e_{\lambda})\omega \varphi \, \mathrm{d}x \right|=0.
    \end{split}
    \label{lim5}
\end{equation}
To this end, we first note
$$
De_{\lambda}  \rightharpoonup D e  \qquad
\textrm{weakly in } L^{\varphi}_{\omega}(\Omega; \mr^n).
$$
Then we discover that there is some $\overline{e}\in
W^{1,q}(\Omega)$ such that
$$
\begin{aligned}
    e_{\lambda} &\rightharpoonup \overline{e} &&\textrm{weakly in } W^{1,q}(\Omega),\\
    D e_{\lambda} &\rightharpoonup D \overline{e} &&\textrm{weakly in } L^{\varphi}_{\omega}(\Omega; \mr^n).
\end{aligned}
$$
Due to the uniqueness of the weak limit, it is enough to check
that $\overline{e}=e$. We have
$$
\begin{aligned}
    \|\overline{e}-e\|_1 &= \lim_{\lambda \to \infty}\int_{\Omega}|e_{\lambda} - e|\, \mathrm{d}x = \lim_{\lambda \to \infty} \lim_{k\to \infty}\int_{\Omega}|e_{k,\lambda} - e_k|\, \mathrm{d}x\\
    &=\lim_{\lambda \to \infty} \lim_{k\to \infty}\int_{\Omega\cap\{M(\nabla e_k)>\lambda\}}|e_{k,\lambda} - e_k|\, \mathrm{d}x \\
    &\le \lim_{\lambda \to \infty} \lim_{k\to \infty} \|e_{k,\lambda} - e_k\|_q |\Omega\cap\{M(D e_k)>\lambda\}|^{\frac{1}{q'}}\le C\lim_{\lambda \to \infty}\lambda^{-\frac{1}{q'}}=0
\end{aligned}
$$
and the equation \eqref{lim5} follows.
 This completes the proof.
\end{proof}

Now we can prove the existence of very weak solution as an
application of div-curl lemma.

\begin{proof}[Proof of Theorem \ref{mainthm}]
We take $\delta_{0} = \frac{\delta_{*}}{2}$ and approximate a given $\mathbf{f} \in L^{\varphi^{1-\delta_{0}}}(\Omega;\mr^n)$ by $$\mathbf{f}_{k}:=
\begin{cases}
    \min \{k,|\mathbf{f}|\} \frac{\mathbf{f}}{|\mathbf{f}|}  &  \textrm{if}\  \mathbf{f}_{k} \neq \mathbf{0}  \\
    0 & \textrm{if}\  \mathbf{f}_{k} = \mathbf{0},
\end{cases} $$ where $\delta_{*}$ is given in Theorem \ref{mainthm2}. Since $\mathbf{f}_{k} \in L^{\infty}(\Omega;\mr^n)$, we can apply standard monotone operator theory in \cite[chapter 2]{Sh} to find the weak solution $u_{k} \in W_{0}^{1, \varphi}(\Omega)$ to the equation
\begin{equation}\label{mainproof21}
  \notag   \begin{cases}
        \mathrm{div\,}  \mathcal{A}(x, Du_{k})  =\mathrm{div\,} \left( \frac{\varphi'(|\mathbf{f}_{k}|)}{|\mathbf{f}_{k}|} \mathbf{f}_{k} \right)  & \textrm{in}\ \Omega  \\
        u=0 & \textrm{on}\ \partial\Omega,
    \end{cases}
\end{equation}
 which means that
\begin{equation}\label{weakform}
\int_{\Omega}       \left\langle   \mathcal{A}(x, Du_{k}), D\eta  \right\rangle    \, \mathrm{d}x   =  \int_{\Omega} \left\langle  \frac{\varphi'(|\mathbf{f}_{k}|)}{|\mathbf{f}_{k}|} \mathbf{f}_{k} ,   D\eta \,  \right\rangle \mathrm{d}x , \quad \forall \eta \in W^{1, \varphi}(\Omega).
\end{equation}
 Then Theorem \ref{mainthm2} along with Lemma \ref{extrapolation} gives
\begin{equation}\label{mainproof22}
     \int_{\Omega}   \varphi(|Du_{k}|)^{q} \, \mathrm{d}x   \leq c  \int_{\Omega}  \varphi(|\mathbf{f}_{k}|)^{q} \, \mathrm{d}x \leq c  \int_{\Omega}  \varphi(|\mathbf{f}|)^{q} \, \mathrm{d}x
 \end{equation}
for every $q \in (1-\delta_{*}, 1+ \delta_{*})$.
Especially, fix $q = 1-\frac{\delta_{*}}{2}$ and take
$$
  \omega := \frac{1}{  \mathcal{M} ( \chi_{\Omega} \varphi (|\mathbf{f} |)^{1-\delta_{*}})^{\frac{\delta_{*}}{2(1-\delta_{*})}} }.   
$$ 
Note that Lemma \ref{maximucken} gives $\omega \in A_{\frac{1}{1-\delta_{*}}} \cap RH_{1+\frac{1}{\delta_{*}}} $. Applying Theorem \ref{mainthm2}, we obtain
\begin{eqnarray}\label{mainproof23}
    \notag &&\int_{\Omega}   \frac{\varphi(|Du_{k}|)}{  \mathcal{M} ( \chi_{\Omega} \varphi (|\mathbf{f} |)^{1-\delta_{*}})^{\frac{\delta_{*}}{2(1-\delta_{*})}} } \, \mathrm{d}x   \leq c  \int_{\Omega}  \frac{\varphi(|\mathbf{f}_{k}|)}{  \mathcal{M} ( \chi_{\Omega} \varphi (|\mathbf{f} |)^{1-\delta_{*}})^{\frac{\delta_{*}}{2(1-\delta_{*})}} } \, \mathrm{d}x  \\
    && \hspace{5mm}\leq c  \int_{\Omega}   \mathcal{M} ( \chi_{\Omega} \varphi (|\mathbf{f} |)^{1-\delta_{*}})^{\frac{1-\delta_{*}/2}{1-\delta_{*}}} \, \mathrm{d}x  \leq c  \int_{\Omega}   \varphi (|\mathbf{f}|)^{1-\delta_{*}/2} \, \mathrm{d}x.
\end{eqnarray}
Therefore, passing to a subsequence, we get
\begin{eqnarray}\label{mainproof24}
     u_{k} \rightharpoonup u \quad \hspace{3mm}     &&\textrm{weakly in}  \, \, W^{1, \varphi^{1-\frac{\delta_{*}}{2}}}(\Omega) \\
 Du_{k} \rightharpoonup Du \quad    &&\textrm{weakly in}  \, \, L_{\omega}^{ \varphi} \cap  L^{\varphi^{1-\frac{\delta_{*}}{2}}} (\Omega; \mr^n)  \\
 \mathcal{A}(x,Du_{k})  \rightharpoonup  \overline{\mathcal{A}}    \hspace{5.5mm}   &&\textrm{weakly in}  \, \, L_{\omega}^{ \varphi^{*}} (\Omega; \mr^n)\label{Aweakconv}
\end{eqnarray}
for some functions $u \in W_{0}^{1, \varphi^{1-\frac{\delta_{*}}{2}}}(\Omega)$ and
 $ \overline{\mathcal{A}} \in L_{\omega}^{ \varphi^{*}} (\Omega; \mr^n)  $.
 Then by \eqref{mainproof22} and  \eqref{mainproof23}, there holds an estimate
\begin{eqnarray}\label{mainproof25}
    \notag  \int_{\Omega}   \varphi(|Du|)^{1-\frac{\delta_{*}}{2}} \, \mathrm{d}x  + \int_{\Omega}   \varphi(|Du|) \omega \, \mathrm{d}x   \leq c  \int_{\Omega}  \varphi(|\mathbf{f}|)^{1-\frac{\delta_{*}}{2}} \, \mathrm{d}x.
\end{eqnarray}
 As $\mathbf{f}_{k}$ converges strongly to $\mathbf{f}$ in
 $L^{\varphi^{1-\frac{\delta_{*}}{2}}} (\Omega; \mr^n)$,
 it directly follows from \eqref{weakform} and  \eqref{Aweakconv}  that
 \begin{eqnarray}\label{mainproof25}
    \notag \int_{\Omega}   \overline{\mathcal{A}} \cdot D\eta \, \mathrm{d}x   =  \int_{\Omega}  \frac{\varphi'(|\mathbf{f}|)}{|\mathbf{f}|} \mathbf{f} \cdot   D\eta \, \mathrm{d}x , \quad \forall \eta \in C_{0}^{\infty}(\Omega).
 \end{eqnarray}
 Then it remains to prove that $\overline{\mathcal{A}} = \mathcal{A}(x,Du) $ in order to conclude that $u$ becomes a distributional solution. Note that
\begin{eqnarray}\label{mainproof26}
   \notag  \int_{\Omega}  | \mathcal{A}(x,Du_{k}) \cdot Du_{k} | \,  \omega \, \mathrm{d}x  \leq  \int_{\Omega}   \varphi(|Du_{k}|) \omega \, \mathrm{d}x   \leq c  \int_{\Omega}  \varphi(|\mathbf{f}|)^{1-\frac{\delta_{*}}{2}} \, \mathrm{d}x.
\end{eqnarray}
To apply divergence curl lemma, Theorem \ref{divcurl}, we take
$a_{k} =  Du_{k}$ and $b_{k} = \mathcal{A}(x,Du_{k})$. Then the
assumption \eqref{divcurlasm2} follows from  \eqref{weakform} while
the assumption \eqref{divcurlasm4} trivially holds since $a_{k}$ is
a gradient. Therefore, we obtain the convergence
\begin{equation}\label{mainproof27}
   \notag  \mathcal{A}(x,Du_{k}) \cdot Du_{k} \, \omega \rightarrow  \overline{\mathcal{A}} \cdot Du \,  \omega
\end{equation}
in distribution sense. Now we apply some variant of the Minty trick.
For any $G\in L^{\varphi}_{\omega}(\Omega; \mr^n)$, we obtain from
the convergence \eqref{mainproof24} that
\begin{equation*}
    (\mathcal{A}(x,Du_{k})-\mathcal{A}(x,G)) \cdot (D u^k-G) \,\omega \rightharpoonup
    (\overline{\mathcal{A}}-\mathcal{A}(x,G)) \cdot (D u-G)\, \omega\quad
\end{equation*}
in distribution sense. Then the monotonicity condition \eqref{moncondi} implies that
\begin{equation}\label{mainproof28}
    \int_{\Omega}(\overline{\mathcal{A}}-\mathcal{A}(x,G)) \cdot (D u-G)\, \eta \omega \, \mathrm{d}x \ge 0
\end{equation}
for any smooth function $\eta \geq 0$. Hence, setting $G:=D u
-\delta H$ with $H\in L^{\infty}(\Omega; \mr^{n})$ to be selected in
a few lines and then dividing \eqref{mainproof28} by $\delta$, we
have
\begin{equation*}
    \int_{\Omega}(\overline{\mathcal{A}}-\mathcal{A}(x,D u - \delta H)) \cdot H\, \omega  \, \mathrm{d}x \ge 0.
\end{equation*}
Finally, recalling the continuity assumption of $\mathcal{A}$ with
respect to $\xi$ variable and using the dominated convergence
theorem, we let  $\delta \to 0$  to discover that
\begin{equation*}
    \int_{\Omega}(\overline{\mathcal{A}}-\mathcal{A}(x,D u)) \cdot H\, \omega \, \mathrm{d}x \ge 0 \qquad \textrm{for all } H\in L^{\infty}(\Omega;\mr^{n}).
\end{equation*}
 Since $\omega$ is strictly positive almost everywhere in $\Omega$,
the relation  $\overline{\mathcal{A}} = \mathcal{A}(x,Du) $ easily
follows by choosing
$$
H:=-\frac{\overline{\mathcal{A}}-\mathcal{A}(x,D u)}{1+|{\overline{\mathcal{A}}-\mathcal{A}(x,D u)}|}.
$$
Therefore, we conclude that $u$ is a very weak solution to
\eqref{maineq3}.
\end{proof}

\bibliographystyle{amsplain}

\end{document}